\documentclass[11pt]{article}
\usepackage{hyperref}
\usepackage{epsfig}

\usepackage{amsmath,amssymb,graphicx}
\usepackage{mathtools}

\usepackage[nosort]{cite}
\usepackage{multicol,color,longtable}

\definecolor{darkred}{rgb}{0.65,0.15,0}
\hypersetup{pdfborder={0 0 0},colorlinks=true,urlcolor=blue,citecolor=blue,linkcolor=darkred,linktocpage=true}

\newcommand{\eprint}[1]{{\href{http://arxiv.org/abs/#1}{[\texttt{#1}]}}}
\newcommand{\eprintN}[1]{{\href{http://arxiv.org/abs/#1}{[\texttt{#1 [hep-th]}]}}}
\newcommand{\eprintQC}[1]{{\href{http://arxiv.org/abs/#1}{[\texttt{#1 [gr-qc]}]}}}

\newcommand{\eprintGR}[1]{{\href{http://arxiv.org/abs/#1}{[\texttt{#1 [math.GR]}]}}}
\newcommand{\eprintRT}[1]{{\href{http://arxiv.org/abs/#1}{[\texttt{#1 [math.RT]}]}}}
\newcommand{\eprintRA}[1]{{\href{http://arxiv.org/abs/#1}{[\texttt{#1 [math.RA]}]}}}
\newcommand{\doi}[2]{{\hypersetup{urlcolor=darkred}\href{http://dx.doi.org/#2}{#1}\hypersetup{urlcolor=blue}}}

\newcommand{\nn}{\nonumber}

\newcommand{\nats}{\mathbb{N}}
\newcommand{\ints}{\mathbb{Z}}
\newcommand{\reals}{\mathbb{R}}
\newcommand{\cx}{\mathbb{C}}
\newcommand{\rats}{\mathbb{Q}}

\newcommand{\cO}{\mathcal{O}}
\newcommand{\cH}{\mathcal{H}}
\newcommand{\cF}{\mathcal{F}}

\newcommand{\cP}{\mathcal{P}}
\newcommand{\cD}{\mathcal{D}}

\newcommand{\lb}{\left[}
\newcommand{\rb}{\right]}
\newcommand{\lp}{\left(}
\newcommand{\rp}{\right)}

\newcommand{\la}{\langle}
\newcommand{\ra}{\rangle}

\newcommand{\ad}{\mathrm{ad}}
\newcommand{\UHP}{\mathbb{H}}

\newcommand{\bi}{{\mathbf{i}}}

\renewcommand{\Re}{\mathrm{Re}}
\newcommand{\sgn}{\mathrm{sgn}}
\newcommand{\bm}{\begin{bmatrix}}
\newcommand{\ebm}{\end{bmatrix}}

\newcommand{\so}{{\mf{so}(2,1)}}
\newcommand{\SO}{{\mf{so}(3)}}

\newcommand{\ba}{{\bar{a}}}

\newcommand{\mf}[1]{{\mathfrak{#1}}}

\newcommand{\beq}{\begin{equation}}
\newcommand{\eeq}{\end{equation}}

\newcommand{\La}{\Lambda}
\newcommand{\al}{\alpha}
\newcommand{\be}{\beta}
\newcommand{\vp}{\varphi}
\newcommand{\N}{|\hspace{-0.5mm}|}
\newcommand{\Herm}{Hermitian }

\numberwithin{equation}{section}

\begin{document}

\mbox{ }\\[10mm]

\begin{center}
{\LARGE \bf \sc Decompositions of hyperbolic Kac-Moody\\[4mm] algebras with respect to imaginary root groups
}\\[18mm]
\end{center}

\begin{center}
Alex J. Feingold\footnotemark[1], Axel Kleinschmidt\footnotemark[2] and
Hermann Nicolai\footnotemark[2]\\[3mm]
\footnotemark[1]{\sl  Department of Mathematics and Statistics, The State University of New York\\ 
 Binghamton, New York 13902--6000, U.S.A.}\\[3mm]
\footnotemark[2]{\sl  Max-Planck-Institut f\"ur Gravitationsphysik,
 Albert-Einstein-Institut \\
 Am M\"uhlenberg 1, D-14476 Potsdam, Germany} \\[7mm]
 \end{center}

\vspace{20mm}

\begin{center} 
\hrule

\vspace{6mm}

\begin{tabular}{p{14cm}}
{\small%
We propose a novel way to define imaginary root subgroups 
associated with (timelike) imaginary roots of hyperbolic Kac--Moody algebras. Using in an essential way the 
theory of unitary irreducible representation of covers of the group $SO(2,1)$, these imaginary root subgroups act on the complex Kac--Moody algebra 
viewed as a Hilbert space. We illustrate our new view on Kac--Moody groups by considering the example of a 
rank-two hyperbolic algebra that is related to the Fibonacci numbers. We also point out some open issues
and new avenues for further research, and briefly discuss the potential relevance of the present results
for physics and current attempts at unification.
}
\end{tabular}
\vspace{5mm}
\hrule
\end{center}

\thispagestyle{empty}

\newpage
\setcounter{page}{1}

\setcounter{tocdepth}{2}
\tableofcontents

\vspace{3mm}
\hrule

\section{Introduction}

The general theory of Kac--Moody (KM) Lie algebras \cite{Kac,Wan} has been recognized as a beautiful and 
natural generalization of the theory of finite-dimensional semi-simple Lie algebras over the complex numbers. 
Many important applications have been found for affine KM Lie algebras, which are characterized
by positive semi-definite Cartan matrices, and their corresponding groups. This statement applies
especially in physics, where the theory has found prominent applications 
in string theory and two-dimensional conformal field theory, but also in the context of axisymmetric stationary 
solutions of Einstein's equations~\cite{Breitenlohner:1986um} and its generalizations. However, the situation with 
indefinite, and more specifically hyperbolic KM algebras, which come with {\em indefinite}
Cartan matrices, is entirely different. These algebras are much more poorly 
understood on the mathematical side, which is mainly due to the presence of imaginary
(time-like) roots whose associated root spaces exhibit exponential growth. As for physical
applications there are tantalizing hints of their possible relevance to understanding the 
physics of the Big Bang in a quantum cosmological context~\cite{Damour:2001sa,Damour:2002cu,Kleinschmidt:2009cv} but it remains unclear how these symmetries 
are to be properly implemented  and interpreted in a physical context. This is even more true with regard to the 
associated KM {\em groups}, where again the  main problem resides with imaginary root space elements
and their exponentiation. 

In this paper we are interested in the hyperbolic KM Lie algebras and groups. So far, KM groups have been mostly defined 
by real root groups generated by exponentials of real root spaces whose adjoint action on the algebra is locally nilpotent.  
The only known definitions of ``imaginary" root groups require completions of the algebra which allow infinite sums in
only one direction (say for positive roots)~\cite{Kumar:2002,Marquis}. In this paper we study an alternative approach which works for one imaginary
root at a time, and uses the extensive theory of unitary irreducible representations of (covers of) $SO(2,1)$. 
This approach has advantages and disadvantages, but we hope that the approach studied here can shed some new 
light on these remarkable algebraic structures which continue to challenge mathematicians and physicists. 

We also believe that potential applications of our results to physics, unification and M-theory are
very interesting. While there is now plenty of evidence that indefinite
KM algebras {\em are} relevant in this context, we have very few tools for dealing with 
them, especially when it comes to the KM {\em groups} obtained by exponentiation 
of the corresponding KM algebras. Even for the KM Lie algebra a physical 
interpretation is so far established only for a finite subset of the real root 
generators and some very specific null roots associated to the elements of the spin 
connection~\cite{Damour:2002cu,Houart:2011sk}.\footnote{See however~\cite{Brown:2004jb} for a 
discussion of some aspects of timelike imaginary roots, and \cite{DN} for partial
evidence associating imaginary roots to higher order corrections in M theory.}  
Likewise the duality symmetries  discussed so far only 
concern finite-dimensional regular subalgebras and their associated 
low level degrees of freedom. By contrast, the $SO(2,1)$ groups 
exhibited here reach `infinitely far' into the space of imaginary root generators,
beyond the low level elements for which a physical
interpretation has been found. If a way could be found to imbue these 
groups with a physical meaning this would open entirely new windows on
string unification, for instance providing new tools to study
higher order corrections beyond perturbation theory.
One important aspect of all proposals including KM symmetries in M-theory is the use of symmetric spaces based on the KM group and the usual physics approach is to include exponentials with all positive root generators, including imaginary root generators. A better understanding of such exponentials was one of the key motivations for this paper.

In the theory of finite-dimensional semi-simple Lie algebras over the complex numbers, 
a great achievement was the Cartan--Killing classification of the simple Lie algebras in terms of an 
integral $n\times n$ Cartan matrix, $A = [a_{ij}]$, which captures the geometry of the root system. Dynkin diagrams are very useful graphs 
which carry the same information as the Cartan matrix, but display it in a clearer way. 
Serre's theorem gives generators, $\{e_i, f_i, h_i\mid 1\leq i\leq n\}$, and relations 
\begin{align}\label{Serre}
&[h_j,e_i] = a_{ij} e_i,\ \ [h_j,f_i] = -a_{ij} f_i,\ \ [e_i,f_j] = \delta_{ij} h_i,\ \ [h_i,h_j] = 0,\ \ \hbox{for } 1\leq i,j \leq n,\\ \nn
&(ad_{e_i})^{1-a_{ij}}(e_j) = 0 =  (ad_{f_i})^{1-a_{ij}}(f_j) \ \ \hbox{for } 1\leq i\neq j \leq n,
\end{align}
for finite-dimensional semi-simple Lie algebras from that Cartan matrix.
Starting from a generalized Cartan matrix, $A = [a_{ij}]$, Kac \cite{Kac:1968} and Moody \cite{Moody:1968} independently in 1968 
defined a class of infinite-dimensional Lie algebras over $\cx$ by generators and Serre relations (see \cite{Gabber-Kac:1981}). 
Most of the results and applications of KM algebras 
have been for the affine KM algebras because they can be described explicitly as a central extension of a loop algebra of a finite-dimensional Lie algebra, 
\beq\label{affineKM}
\hat {\mf{g}} = \mf{g}\otimes \cx[t,t^{-1}] \oplus \cx c \oplus \cx d
\eeq
(the ``untwisted case'') where $\mf{g}$ is a finite-dimensional semi-simple Lie algebra over $\cx$, $c$ is central and $d$ is a derivation acting on the ring
$\cx[t,t^{-1}]$ needed to extend the Cartan subalgebra because the affine Cartan matrix $A = [a_{ij}]$ has $\det(A) = 0$. 
Lie brackets for these are explicitly given, in contrast with the indefinite KM algebras where the definition only gives a generators and 
relations description. 

Among the indefinite KM algebras, the class of hyperbolic type has received the most attention, including some 
applications in theoretical physics to supergravity. The representation theories of the affine and hyperbolic types are also in stark contrast
mainly because each affine KM algebra contains an infinite-dimensional Heisenberg Lie subalgebra,
\beq\label{Heisenberg}
\hat {\mf{h}} = \mf{h}\otimes \cx[t,t^{-1}] \oplus \cx c
\eeq
where $\mf{h}$ is the abelian Cartan subalgebra of $\mf{g}$. The Fock-space representation of $\hat {\mf{h}}$ on a space of polynomials in 
infinitely many variables as multiplication and partial differentiation operators plays a vital role in the vertex operator representations of 
$\hat {\mf{g}}$. The rank $2$ hyperbolic KM algebras do not contain any Heisenberg Lie subalgebra. For higher rank hyperbolic KM algebras
which contain an affine KM subalgebra, one can decompose the hyperbolic algebra with respect to its affine subalgebra or some other kind of 
subalgebra whose representations can be understood. 

There have been several choices studied for such a subalgebra in a hyperbolic $\mf{g} = \mf{g}(A)$: 

(1) A finite type KM algebra coming from a subset of the generators (a Dynkin sub-diagram), 

(2) An affine type KM algebra coming from a subset of the generators, 

(3) A subalgebra of fixed points under an automorphism of $\mf{g}$,

(4) A subalgebra which is not obvious, e.g., not just from a Dynkin sub-diagram.

Option (1) has been used, for example, to study the hyperbolic algebra known as $E_{10}$ by decomposing it with respect to a  
finite type $A_9$ subalgebra~\cite{Damour:2002cu,Nicolai:2003fw}. 
Similar decompositions have been performed with respect to the $D_9$ and $A_8 \oplus A_1$ subalgebras in~\cite{Kleinschmidt:2004dy,Kleinschmidt:2004rg}.
In physical applications a real Lie algebra is preferred, usually the split real form, $\mf{g}_\reals = \mf{g}(A)_\reals$,
which is just the real span of the generators and their Lie brackets. The split real form can also be understood under option (3) as the fixed points 
of the conjugate linear involutive automorphism that fixes the Chevalley generators, $\{e_i, f_i, h_i\mid 1\leq i\leq n\}$.  

Option (2) has been used, for example in \cite{Feingold-Frenkel:1983}, to study a particular rank 3 hyperbolic, $\cF$, also called 
$AE3$, which has an affine subalgebra of type $A_1^{(1)}$, the simplest example of an affine KM algebra whose representation theory is well developed. 

Option (3) includes the split real form mentioned above, as well as the ``compact'' real form, $K(\mf{g})$, which is a real Lie subalgebra of fixed points 
under the Cartan--Chevalley involution $\omega(e_i) = -f_i$, $\omega(f_i) = -e_i$ and $\omega(h_i) = -h_i$ on the complex KM algebra. 
The intersection of the split and the compact real form is of interest to physicists~\cite{Nicolai:2004nv,deBuyl:2005zy,Damour:2005zs,deBuyl:2005sch,Damour:2006xu,Kleinschmidt:2021agj}, 
who have studied finite-dimensional representations of the infinite-dimensional 
involutive subalgebra $K(\mf{g}_\reals)$ generated by $\{k_i = e_i - f_i \mid 1\leq i\leq n\}$ and satisfying the Berman relations \cite{Berman:1989}. 
When a finite type algebra has a Dynkin diagram with an automorphism (symmetry), twisted affine KM algebras result from the fixed point subalgebra of $\hat {\mf{g}}$. 

Option (4) can be applied using the results of \cite{Feingold-Nicolai:2004} on subalgebras of hyperbolic KM algebras. They found inside $\cF$ all the 
rank 2 hyperbolics whose Cartan matrix is symmetric. The simplest example is the rank 2 ``Fibonacci'' hyperbolic $Fib$ \cite{Feingold:1980} whose $2\times 2$ Cartan 
matrix has $a_{12} = -3 = a_{21}$. A study was made in \cite{Penta:2016} of the decomposition of $\cF$ with respect to $Fib$ that showed some 
interesting $Fib$-modules occur, including some integrable modules which are neither highest nor lowest weight modules, and not the adjoint 
module. We will later use $Fib$ as one of the simplest examples of a hyperbolic KM algebra to illustrate how it and two of its irreducible 
highest weight representations (see Figures~\ref{fig:Vrho} and~\ref{fig:Vlam1}) might be decomposed in a new way using a three-dimensional {\it imaginary} subalgebra determined by a choice of 
an imaginary root vector in some imaginary root space of a hyperbolic KM algebra. For comparison we will also discuss how a choice of a real 
root vector in a real root space gives a decomposition into finite-dimensional $\mf{sl}(2,\cx)$-modules. The use of ``real'' versus ``imaginary'' for 
kinds of roots should not be confused with the choice of field $\reals$ versus $\cx$ for the scalars of the Lie algebra and its representations. 

Included in option (1) is the obvious choice of an $\mf{sl}(2,\cx)$ subalgebra corresponding to a simple root, $\al_i$, for a fixed $i$, that is, the 
subalgebra $\mf{g}_i = \mf{sl}(2,\cx)_i$ with basis $\{e_i, f_i, h_i\}$. The Serre relations defining $\mf{g}$ imply that it decomposes with respect to $\mf{g}_i$
into an infinite number of finite-dimensional $\mf{g}_i$-modules. We could have taken any real root, $\al$, whose root space $\mf{g}_\al$ must be 
one-dimensional with basis vector $E(\al)$, and found an opposite root vector $F(\al)$ in $\mf{g}_{-\al}$, such that with $H(\al) = [E(\al),F(\al)]$
a subalgebra $\mf{sl}(2)_\al$ is defined. But since any real root is by definition in the Weyl group orbit of the simple roots, it is sufficient to study
just the decompositions with respect to the subalgebras $\mf{sl}(2)_i$. In a later section we will discuss in some detail how this decomposition works
for the rank 2 hyperbolic algebra $Fib$, see section~\ref{sec:Fibadj}. See Figure~\ref{fig:Fibpos} for a graphical display of some positive roots of $Fib$ along with their root multiplicities.

Our construction raises several interesting questions and opens new avenues for further
research. One of them concerns the issue of `combining' different $SO(2,1)_\alpha$ groups and their interplay
for {\em different} imaginary roots $\al$. Unlike for real root subgroups, there are no (Steinberg-type) relations that could be
exploited towards the evaluation of products of elements of different imaginary root groups due to the lack of local nilpotency.\footnote{We note that, according to the results of~\cite{Marquis:2021}, generators associated with root spaces of different positive imaginary roots generate a free Lie algebra under some mild assumptions.}
Although each action would involve distinct Hilbert spaces, the repeated action of such operations 
is well-defined, because the unitary action guarantees that norms are preserved by the repeated
group action.  

What we would like to stress here that for the definition of imaginary root subgroups {\em some}
notion of completion is definitely required. The one we employ here relies on the Hilbert space completion of 
the vector space of the KM algebra viewed as an $SO(2,1)_\alpha$ module. To what extent the KM algebra 
structure is compatible with this completion remains an open question, 
but the results of~\cite{Kleinschmidt:2021agj} suggest the thus completed space may no longer be a KM algebra,
in the sense that the commutator of two elements of the completion is
no longer an element of the Hilbert space.
This could mean that the Hilbert space norm used here which is induced by
the standard bilinear form may not be the appropriate tool to define a completion of the KM algebra,
and that one may have to resort to different notions of completion (we note that there are 
many topologies on infinite dimensional vector spaces that might be used here).
This is also exemplified by comparing  between Kac--Moody commutators and tensor products in the general theory of unitary 
representations of $SO(2,1)$ as described for instance in \cite{Pukanszky:1961,Repka}. More specifically,
take the two principal series that arise in the adjoint of $Fib$ derived in section~\ref{sec:Fibadj}. They are both unitary but their commutator  contains, among other things, the {\em non-unitary} 
adjoint of $\mathfrak{so}(2,1)$. 
This is in tension with the tensor product results 
given in~\cite{Pukanszky:1961,Repka}, according to which the product of two unitary
representations is again unitary. 
A possible explanation of this tension is that the norm of a commutator in the Kac--Moody algebra is not equal to the product of the norms of its two elements which underlies the tensor product construction in the general theory of \cite{Pukanszky:1961,Repka}.

\subsubsection*{Acknowledgements}
We are grateful to Lisa Carbone, Thibault Damour, Walter Freyn, Benedikt K\"onig, Robin Lautenbacher and Timoth\'ee Marquis for discussions and the referees for useful comments. AF gratefully acknowledges 
support from his department and from the Max Planck Institute for Gravitational Physics during several visits related to this work. This work was supported in part by the European 
Research Council (ERC) under the European Union's Horizon 2020 research and innovation programme (grant agreement No 740209).

\section{Decompositions of hyperbolic Kac--Moody algebras}

In this section we will discuss the main idea of the paper, how the choice of an imaginary root vector (multi-bracket) gives an imaginary three-dimensional
subalgebra whose split real form is isomorphic to $\so$. We use the representation theory of $\mf{sl}(2,\reals)\cong \so$ on well-known series of 
unitary modules to decompose any hyperbolic Kac--Moody algebra or representation in such a way that the action of the group $SO(2,1)$ (or its covers) is given explicitly on 
a Hilbert space completion of each irreducible summand. 
This approach defines imaginary root groups in a different way from other methods that use a completion of the Kac--Moody group in 
only one ``direction'', see e.g.~\cite{Marquis} and other references on page 268 of that book or~\cite{Kumar:2002}.

\subsection{Kac--Moody algebras and involutions}
\label{sec:KMAs}

Let $\mf{g} = \mf{g}(A)$ be a Kac--Moody (KM) algebra  over $\cx$ with a non-degenerate symmetric $r\times r$ 
Cartan matrix $A = [a_{ij}]$. The generalization of this work to symmetrizable
Cartan matrices should be straightforward. Our main focus will be on hyperbolic KM algebras where it is well-known that 
$A$ is Lorentzian with signature $(1,r-1)$ and the maximal rank $r$ is $10$ \cite{Kac}, but these introductory remarks are valid in greater generality.

We recall that $\mf{g}$ has a root space decomposition
\begin{align}
\mf{g} = \mf{h} \oplus \bigoplus_{\alpha\in\Delta} \mf{g}_\alpha
\end{align}
where $\mf{h}$ is the $r$-dimensional Cartan subalgebra (CSA) that acts semi-simply by the adjoint action $\ad_h(x) = [h,x]$ 
on $\mf{g}$, and $\mf{g}_\alpha$  denotes the eigenspaces under this action with eigenvalue given by roots $\alpha\in\Delta$ of the form 
$\alpha: \mf{h} \to \mathbb{C}$ such that the eigenspace is non-trivial. Roots $\alpha$ are divided into real roots (characterized by positive norm squared) 
and imaginary roots. The latter can be further subdivided into lightlike (with vanishing norm squared) and timelike roots (with negative norm squared).

The Cartan--Chevalley involution is the $\cx$-antilinear automorphism of $\mf{g}$ ($\omega(zx) = {\bar z}\omega(x)$ for $z\in\cx$ and $x\in\mf{g}$) defined by
\begin{equation}
\omega(e_i) = - f_i \;,\quad \omega(f_i) = - e_i \;,\quad \omega(h_i) = - h_i 
\end{equation}
and extended to the whole KM algebra  by means of $\omega([x,y]) = [\omega(x),\omega(y)]$.
In particular, for any multi-bracket we have
\beq
\omega\big( e_{i_1\cdots i_n}\big) = (-1)^n  f_{i_1\cdots i_n}
\eeq
where we use the notation $e_{i_1\ldots i_n} \coloneqq [e_{i_1},[e_{i_2},\dots , [e_{i_{n-1}} ,e_{i_n}]...]]$ and similarly for  
$f_{i_1\cdots i_n}$. The standard bilinear form is defined by
\beq\label{stdbilform}
\langle e_i | f_j \rangle = \delta_{ij} \;\;,\quad \langle h_i | h_j \rangle = a_{ij}
\eeq
and $\langle [x,y]| z\rangle = \langle x | [y,z]\rangle$.
Then Theorem 11.7 of \cite{Kac} shows that the {\em \Herm  form} (complex-conjugate linear in the second argument)
\beq\label{HF}
(x,y) \coloneqq - \big\langle x|\omega(y)\big\rangle
\eeq
is positive definite on the whole (complex) KM algebra except on its Cartan subalgebra, where it has precisely
one negative eigenvalue. For any operator $\cO$ on $\mf{g}$, its \Herm  conjugate $\cO^\dagger$ is defined by $(\cO(x),y) = (x,\cO^\dagger(y))$
for any $x,y\in\mf{g}$. 

With respect to this \Herm form, for any element $z\in \mf{g}$, the adjoint operator $\ad_z$, defined by $\ad_z(x) = [z,x]$ for any $x\in\mf{g}$,
satisfies 
\beq
(\ad_z)^\dagger = \ad_{-\omega(z)}. 
\eeq
To see this we check
\begin{align}
\big(\ad_z(x),y\big) = \big([z,x],y\big) &= - \big\langle [z,x]\,|\,\omega(y)\big\rangle = + \big\langle x\,|\,[z,\omega(y)]\big\rangle 
\nn\\
&= \big\langle x\,\big|\, \omega \lb\omega(z),y\rb\big\rangle  = \big(x\,,[-\omega(z),y] \big) = \big(x\,,\ad_{-\omega(z)}(y) \big).
\end{align}
In particular, $\ad_z$ is self-conjugate if and only if $z = -\omega(z)$.

\subsection{\texorpdfstring{Subalgebras $\so_\al$ associated with  roots $\al$}{Subalgebras so(2,1) associated with  roots alpha}}
\label{sec:so21s}

The following works for any indefinite KM algebra $\mathfrak{g}$ whose Cartan matrix has indefinite signature, not only for hyperbolic algebras.

Let $\al = \al_{i_1} + \cdots + \al_{i_n}$ be a positive root belonging to some multi-commutator $e_{i_1\cdots i_n}$, 
where each $\al_{i_j}$ is a simple root and the $e_{i_j}$ are the Chevalley generators. The order of indices
is significant if $\alpha$ is not a real root. Define 
\beq
\label{eq:EFdef}
E(\al) \coloneqq\,  e_{i_1\cdots i_n} \,\in\, \mf{g}_\al \;\;, \quad 
F(\al) \coloneqq - \omega(E(\al)) = - (-1)^n f_{i_1\cdots i_n} \,\in\, \mf{g}_{-\al}
\eeq
so that 
\beq 
N\coloneqq \langle E(\al)| F(\al)\rangle = (E(\al),E(\al)) > 0
\eeq 
In principle we should use the multi-index label $(i_1,\ldots, i_n)$ instead of just $\alpha$ to distinguish the 
independent elements of the root space $\mf{g}_\al$, but we suppress this for simplicity of notation.\footnote{We could also take linear combinations of {\em different} elements of $\mf{g}_\al$, but that would not affect the main argument.} 
The fact that $\alpha$ is the root of $E(\alpha)$ means that
\beq
\big[h_j, E(\alpha) \big] = \al(h_j) E(\alpha)\ .
\eeq
Writing $\al=\sum_{j=1}^r n^j\al_j$ and defining $H(\alpha) = \sum_{j=1}^r n^j h_j$ we also have
\begin{align}\label{EFH}
\big[ H(\al)\,,\, E(\al) \big] &= \al^2 E(\al) \nonumber\\
\big[ H(\al)\,,\, F(\al) \big] &= - \al^2 F(\al)  \nonumber\\
\big[ E(\al)\,,\, F(\al)\big] &= N H(\al)\,.
\end{align}
with $\alpha^2 = \sum n^i a_{ij} n^j$ and where the last equation uses the invariance of  the standard bilinear form $\langle [H(\al),E(\al)]| F(\al)\rangle = \langle H(\al)|[E(\al), F(\al)]\rangle$ with the normalization $\langle h_i | h_j\rangle = a_{ij}$.

Now we have to distinguish two cases. When $\alpha$ is a real root
($ \alpha^2 \coloneqq \alpha\cdot\alpha >0$), we define
\beq\label{JJJ1}
J_3 = (\al^2)^{-1} H(\al) \;,\quad 
J^+ = (N\al^2)^{-1/2} E(\al) \;,\quad J^- = (N\al^2)^{-1/2} F(\al) \;\;,\quad
\eeq
The commutation relations are
\begin{align}
\label{SO3}
\lb J^+ , J^- \rb = + J_3\,,\quad \lb J_3, J^{\pm} \rb = \pm J^{\pm} \,.
\end{align}
These are elements of the KM algebra, so we understand them as operators under the adjoint action.
The hermiticity properties of these operators are inherited from the bilinear form, that is, 
with respect to the Cartan--Chevalley involution $\omega$ the generators satisfy
\begin{align}
\omega( J^\pm) =  - J^{\mp}\,,\quad \omega(J_3) =  - J_3
\end{align}
whence we have
\begin{align}\label{Herm}
\lp J^\pm \rp^\dagger = J^\mp\,,\quad \lp J_3\rp^\dagger = J_3 \,.
\end{align}
One easily checks that these elements have positive norm with respect to the \Herm form (\ref{HF}). 
The generators $J^\pm$ and $J_3$  together with 
the commutation relations (\ref{SO3}) and the hermiticity properties (\ref{Herm}) 
therefore represent the real Lie algebra $\SO$. For real roots, the (adjoint) action of 
$J^\pm$  on the KM algebra $\mf{g}(A)$ generates finite-dimensional 
representation spaces because the multiple addition of 
a real root $\al$ to any root $\be$ will satisfy $(\be + k\al)^2 > 2$ for sufficiently large $k$. 
The associated groups obtained by exponentiating these Lie algebra
elements are referred to as {\em real root groups}, where the exponentiation
can be performed over $\reals$ or $\cx$, or any other field of characteristic zero. These real root groups generate the {\em minimal Kac--Moody group}
associated with the Cartan matrix~\cite{Peterson:1983,Kumar:2002,Marquis}.

The second case to be considered concerns imaginary roots, for which $\al^2\leq 0$.
For lightlike imaginary roots $\alpha$, for which $\alpha^2 = 0$,
one obtains a Heisenberg algebra from~\eqref{eq:EFdef} that corresponds to a contraction of $\mf{sl}(2)$.
However, our main interest here is the case of {\em timelike} imaginary roots,
for which $\alpha^2 < 0$. In that case we can define a subalgebra  $\so_\al$ 
of the KM algebra $\mf{g}(A)$ for {\em any} element of a timelike imaginary root space $\mf{g}_\al$. 

Instead of (\ref{JJJ1}), the relevant definition reads now for timelike roots
\beq\label{JJJ}
J_3 = (\al^2)^{-1} H(\al) \;,\quad 
J^+ = (-N\al^2)^{-1/2} E(\al) \;,\quad J^- = (-N\al^2)^{-1/2} F(\al) \,.
\eeq
It is straightforward to see that these operators satisfy the bracket relations of 
an $\so$ Lie algebra, that is, 
\begin{align}
\label{SO2,1}
\lb J^+ , J^- \rb = - J_3\,,\quad \lb J_3, J^{\pm} \rb = \pm J^{\pm} \,,
\end{align}
which differs by a crucial minus sign from (\ref{SO3}) in the first commutator, 
while the hermiticity 
properties (\ref{Herm}) are maintained. The latter point is essential, since
otherwise the minus sign could simply be redefined away, for instance by rescaling 
$J^+ \rightarrow - J^+$, but this redefinition would violate the hermiticity properties
(\ref{Herm}).
The normalization (\ref{JJJ}) implies that for $\al^2 < 0$
\beq\label{Jnorm}
\N J^+ \N^2 = \N J^-\N^2 = - (\al^2)^{-1}  > 0 \;\;,\quad \N J_3\N^2  = (\al^2)^{-1} < 0
\eeq 
so these norms shrink to zero as $\al^2 \rightarrow - \infty$.

The difference between the real Lie algebra $\mf{so}(2,1)_\al$ for timelike roots compared to $\mf{so}(3)$ for real roots becomes apparent 
when writing these algebras in terms of standard Lorentz or rotation algebras as reviewed in appendix~\ref{sec:sl2}. 
The change of basis from the standard basis to~\eqref{SO3} or~\eqref{SO2,1} involves complex coefficients in such a way that the hermiticity 
properties of the algebras are different in unitary representations. This will also be important when considering the implications for the Kac--Moody 
group in section~\ref{sec:KMgp}.
Since the definition of the generators in~\eqref{JJJ} depends on the root $\alpha$,
we will keep this dependence in the notation for the algebra  $\so_\al$.
As a real Lie algebra we have the isomorphism $\mf{so}(2,1)_\alpha \cong \mf{sl}(2,\mathbb{R})$.

Before continuing we note that there is another way to define $\so$ subalgebras
of $\mf{g}(A)$ in the case of hyperbolic algebras that does not make use of timelike imaginary roots, but rather 
appropriate linear combinations of {\em real} roots. Distinguished among these
is the  principal $\so$ subalgebra introduced in \cite{Nicolai:2001ir}. 
Generalizations  of this construction are studied in \cite{Tsu}.
The principal $\so$ subalgebra  can be constructed using the inverse Cartan matrix 
$A^{-1} = [b_{ij}] = [\La_i\cdot\La_j]$, where $\La_i$, $1\leq i\leq r$,  are the fundamental weights. 
The entries of $A^{-1}$ satisfy $b_{ij} \leq 0$ since all fundamental weights are null or 
time-like for hyperbolic KM algebras, so their scalar products are non-positive. 
We recall that we assume the Cartan matrix to be symmetric for simplicity. If we define
\begin{align}
r_i = - \sum_j b_{ij} >0
\end{align}
then the generators
\beq
\label{eq:princso21}
J_3 = - \sum_i r_i h_i = \sum_{i,j} b_{ij} h_i\,,\qquad
J^+ = \sum_i \sqrt{r_i} e_i \,,\qquad 
J^- = \sum_i \sqrt{r_i} f_i 
\eeq
again satisfy the commutation relations (\ref{SO2,1}) and the 
hermiticity properties (\ref{Herm}).

\subsection{\texorpdfstring{Decomposing $\mf{g}(A)$ under the action of $\so_\al$}{Decomposing g(A) under the action of so(2,1)}}

The subalgebra $\so_\al\subset\mf{g}$ can be used to decompose the adjoint representation (or any other representation) of $\mf{g}$ under its action. Since $\so_\al\cong\mf{sl}(2,\mathbb{R})$ we will be dealing with representations of $\mf{sl}(2,\mathbb{R})$. In view of the hermiticity properties~\eqref{Herm} these representations will typically be unitary representations of $\mf{sl}(2,\mathbb{R})$ so we review the relevant infinite-dimensional representation spaces, called principal series and discrete series representations, in appendix~\ref{sec:absrep}. From now on we  take $\mf{g}$ to be a complex hyperbolic KM algebra.

To analyse the decomposition of the adjoint $\mf{g}$ under the algebra  $\so_\al$ generated
by (\ref{EFH}) let us consider an arbitrary imaginary root $\be$ and any element
$E(\be) \in \mf{g}_\be$ of its associated root space; then
\beq\label{J3beta}
\big[ J_3 \,,\, E(\be)\big] \,=\,   \nu \, E(\be)\,,
\quad\text{where}\quad \nu= \frac{\al\cdot\be}{\al^2} \in \mathbb{Q}\,.
\eeq
For $\be$ a positive timelike imaginary root we have $\al\cdot\be <0$ and
therefore the prefactor on the right-hand side is positive. 
In general, the rational number $\nu$ is not an integer. 
While this does not matter much for the representations of the Lie algebra $\mf{sl}(2,\mathbb{R})$, this matters for the group: the exponential operator $e^{2\pi i r J_3}$ is not periodic modulo $2\pi$ if the $J_3$ eigenvalue $\nu$ not an integer.
In other words, the group obtained by exponentiation of $\so_\alpha$ is not $SO(2,1)$ but a covering of it.  We note that the parameter $\nu$ can become arbitrarily small.
There are infinitely many covers since the fundamental group $\pi_1(SO(2,1)) \cong \ints$ and this agrees with the fact that any denominator can occur in $\nu$ as $\alpha$ varies. The most well-known cover is $Spin(2,1)$ corresponding to the double cover (with the metaplectic Weil representation of $SL(2,\mathbb{R})$) but more complicated situations are possible. In physical terms, such 
representations are {\em anyonic representations of covers of} $SO(2,1)$~\cite{Froehlich,Froehlich1,Luescher,Jackiw,Cortes:1994ff,Cortes:1995wa}.\footnote{The fundamental group of $SL(2,\reals)$ is also $\ints$. Bargmann's classification~\cite{Bargmann:1946me} only addresses representations of $SL(2,\reals)$, not of its higher covers.} There are such representations both for the so-called principal series, the discrete series and the complementary series.

Returning to the decomposition of $\mf{g}$ under $\so_\al$ for positive timelike $\al$, we note that, with respect to the bilinear form \eqref{stdbilform}, the orthogonal 
complement in the Cartan subalgebra $\mf{h}$ of $J_3$ from \eqref{JJJ} consists of singlets. In particular, choosing a basis 
of $(r-1)$ CSA generators $H(v_i)$ with $\al\cdot v_i =0$ (with space-like $v_i$), 
we have $[J^\pm, H(v_i)] = 0$; all these states have positive norm because the $v_i$ are spacelike. 

For other representations let us pick any positive root $\be$, and apply $J^\pm$ to
any element $E(\be)\in\mf{g}_\beta$. The successive application of the lowering operator 
$J^-$ will result in a chain of maps
\begin{align}\label{Jminus}
\ldots \overset{J^-}{\longrightarrow} \mf{g}_{\beta} \overset{J^-}{\longrightarrow} \mf{g}_{\beta-\alpha} \overset{J^-}{\longrightarrow} \mf{g}_{\beta-2\alpha} \overset{J^-}{\longrightarrow} \mf{g}_{\beta-3\alpha} \overset{J^-}{\longrightarrow} \ldots
\end{align}
along an infinite string of subspaces 
\beq\label{gbeta}
\bigoplus_{k\in\ints} \, \mf{g}_{\beta + k\alpha} \subset \mf{g}(A)
\eeq
Likewise, the application of $J^+$ moves in the opposite direction:
\begin{align}\label{Jplus}
\ldots \overset{\;\;\; J^+}{\longleftarrow} \mf{g}_{\beta+3\alpha} 
\overset{\;\;\; J^+}{\longleftarrow} \mf{g}_{\beta+ 2 \alpha} \overset{\;\;\; J^+}{\longleftarrow} 
\mf{g}_{\beta+ \alpha} \overset{\;\;\; J^+}{\longleftarrow} \mf{g}_{\beta} 
\overset{\;\;\; J^+}{\longleftarrow} \ldots
\end{align}
For the root string $\{\be + k\al \mid k\in\ints\}$ we must distinguish two main cases:
\begin{itemize}
\item there exists a minimal $k_0  < 0$ such that $\be  + k_0\al$ is not a root;
\item the elements $\be + k\al$ are roots for all $k\in\ints$ (the chain may or may not
         contain a real root)
\end{itemize}
In the first case the chain terminates and all elements of the root spaces along the
chain belong to discrete representations (idem for the negative side). In the second
case we may encounter continuous representations. However, also in that case
there will occur (many!) discrete representations in the subspace~\eqref{gbeta} of $\mf{g}$. This is because the root multiplicities, 
defined as mult$(\beta) = {\rm dim} \,\mf{g}_\beta$, vary with $\beta$, and
increase exponentially with the height of the root. 
Namely, the subspaces $\mf{g}_{\be + k\al}$ are in general of different dimensions,
with multiplicities increasing in the leftward direction for positive roots (as long as
$\be + k\al$ is positive), and likewise in the rightward direction for negative roots. 
This implies that for positive $\be + k\al$ each root space in the descending chain 
(\ref{Jminus}) has a large kernel whose elements are annihilated by the action of $J^-$.
Consequently, for each root space, every element of the kernel is a lowest weight 
vector of a discrete series representation that extends to the left and is generated 
by the successive application of $J^+$, with the value of $m$ given by formula (\ref{J3beta}). Because $\al\cdot\be/\al^2 >0$ for all positive $\be$ the unitarity condition 
(\ref{s}) discussed in the appendix is  satisfied.
This number is in general fractional and thus we are dealing with
an {\em anyonic} discrete series representation of a cover of $SO(2,1)_\al$.
This covering has at most $(-\al^2)$ sheets, but since $(-\al^2)$ can become arbitrarily 
large we eventually reach all positive rationals. 
A natural framework for considering all possibilities of imaginary timelike roots together is therefore the universal cover.

If the chain of roots $\be + k\al$ does not terminate we are left with principal series representations after the elimination of the discrete series representations. 
In the examples we have studied so far, we noticed that no complementary series representations occurred in the decompositions and this is discussed further in Section~\ref{sec:Fibimag2}. In the following, we shall assume that they are always absent for simplicity although their presence would not qualitatively change our structural analysis.

We do have a proof that there can only be a finite number of principal series representations since any principal series has a string of weights $\beta+k\alpha$ with 
$k\in\mathbb{Z}$. This string must therefore intersect the region in $\mathfrak{h}^*$ bounded by two planes orthogonal to the timelike $\alpha$ that are separated by $\alpha$. 
The intersection of this region with $\Delta\cup\{0\}$ only contains finitely many elements and therefore only  finitely many Lie algebra generators can be part of principal series 
representations, showing that these are finite in number.
The determination of the Casimir for a principal series representation must be done `by hand', as there appears to be no general
formula for evaluating the requisite multi-commutators. 
In section~\ref{sec:Fibadj} we will therefore present a few exemplary calculations involving the so-called ``Fibonacci'' algebra, $Fib$,
which is the simplest example of a (strictly) hyperbolic KM algebra.

In total, we obtain therefore the following decomposition of the vector space of the complex KM algebra $\mf{g}$ as modules for $\so_\al$:
\begin{align}\label{so_alpha_decomp_of_g}
\mf{g} = (\so_\al \otimes \mathbb{C}) \oplus \mathbb{C}^{r-1} \oplus \bigoplus_{s\in Pr_\al(\mf{g})} Mult^{pr}_s \mathcal{P}_s
\oplus \bigoplus_{s> 0} Mult^+_s \mathcal{D}^{+}_s   \oplus \bigoplus_{s>0} Mult^-_s \mathcal{D}^{-}_s 
\end{align}
The first term is the complexified adjoint of $\so_\al$, the second term represents the singlets in the CSA, the third term is a sum over the set $Pr_\al(\mf{g})$ of principal 
series representations which occur in this decomposition with multiplicity $Mult^{pr}_s$. 
The last two terms are the infinite number of lowest and highest weight discrete series representations, each 
labelled by parameter $s$, and the number of copies for each value of $s$ given by the multiplicity $Mult^\pm_s$. Since we are dealing with covers of $SO(2,1)$, 
only restriction on the value of $s$ for discrete series representations is being a positive real number. 
For the group $SO(2,1)\cong PSL(2,\mathbb{R})$, we would require a positive integer.  The various multiplicities which occur in this decomposition are computable in
any particular example from a knowledge of the root multiplicities of $\mf{g}$, but we have no general formula for them. 
It follows from the Cartan--Chevalley involution that there is a bijection between the lowest and highest weight representations that arise, so we get $Mult^+_s = Mult^-_s$. 
The representation spaces are normed and carry a unitary action of $\so_\al$, hence they have completions, $\mathcal{\hat P}_s$ and $\mathcal{\hat D}^{\pm}_s$, 
which are complex Hilbert spaces. While this gives a Hilbert space completion of the vector space, $\mf{g}$, for each choice of $\al$, we do not claim that the Lie algebra
structure of $\mf{g}$ extends to the completion. We have used the notation introduced in appendix~\ref{sec:absrep}. 
In section~\ref{sec:KMgp} this will be important since we know that we can define the action of a group (a cover of $SO(2,1)$) on these spaces.

We may also decompose any highest weight representation, $V = V(\Lambda)$, of $\mf{g}$ with respect to $\so_\al$:
\begin{align}\label{so_alpha_decomp_of_V}
V(\Lambda) = \bigoplus_{s\in Discr(V)} Mult^-_s \mathcal{D}^{-}_s 
\end{align}
where $Discr(V)$ is the set of values of parameter $s$ which occur in the decomposition, and $Mult^-_s$ is the multiplicity giving the number of copies which occur
of the irreducible highest weight discrete series representation $\mathcal{D}^{-}_s$. A similar formula holds for lowest weight representations of $\mf{g}$ where only
lowest weight discrete series representations $\mathcal{D}^{+}_s$ occur. In section~\ref{sec:Fibadj} we will give some terms in the decompositions of two highest weight
representations of $Fib$. In order to get representations of covers of $SO(2,1)_\al$, these discrete series representations must be completed with respect to a norm to 
Hilbert spaces, $\mathcal{\hat D}^{\pm}_s$, so there is a Hilbert space completion, ${\hat V}(\Lambda)$, of $V(\Lambda)$. Now the question is whether this completion,
which depends on the choice of $\al$, is a $\mf{g}$-module.

\medskip

To close this section, we give the decomposition of $\mf{g}$ under the principal $\so$ subalgebra~\eqref{eq:princso21} that was already studied
in \cite{Nicolai:2001ir}. There it  was shown to take the following form 
\begin{align}\label{decomp1}
\mf{g} = (\so\otimes \mathbb{C}) \oplus \bigoplus_{i=1}^{{r}-1} \mathcal{P}_{i}
 \oplus \bigoplus_{k\geq 2} Mult^+_k \mathcal{D}_k^{+}\oplus \bigoplus_{k\geq 2} Mult^-_k\mathcal{D}_k^{-}\,,
\end{align}
where there are $r -1$ unitary principal series representations $\mathcal{P}_{i}$ whose $s$ parameter is suppressed in this formula, 
and $Mult^\pm_k$ discrete series representations $ \mathcal{D}_k^{\pm}$ of lowest ($+$) or highest ($-$) weight type and parameter $k$ with $2\leq k\in\ints$. 
By the Cartan--Chevalley involution there is a bijection between the set of lowest weight representations and the set of highest weight representations, 
so $Mult^+_k = Mult^-_k$. Because they appear with \textit{integral} parameters, they lift to representations of $SO(2,1)$. Only single-valued 
representations appear in (\ref{decomp1}), but the decompositions take a more 
complicated form for the $\so_\al$ subalgebras associated with timelike imaginary roots $\al$.

\section{Defining a group action on the Kac--Moody algebra}
\label{sec:KMgp}

A main goal of this paper is to see under what circumstances an exponential
action of the imaginary root subalgebras $\so_\al$ on the full KM algebra can be defined. 
For {\em real root subalgebras} this is always possible, yielding an action of the respective {\em real 
root subgroups}. This can be done with either real or complex coefficients (or even coefficients
taking values in some ring or finite field)  because the relevant representations are always finite-dimensional.
This is the reason why the minimal KM groups are defined to be generated by the real root subgroups. 
However, defining larger KM groups including generators from the imaginary root subalgebras $\so_\al$ is only possible
subject to special restrictions. One option that has been extensively explored in the literature (see \cite{Kumar:2002,Marquis} and references therein)
is to consider only exponentials of {\em positive} imaginary root vectors. 
Here, we wish
to explore an alternative option that allows us to exponentiate the action of any imaginary subalgebra 
$\so_\al$ to get an imaginary root group action on completions of $\mf{g}$ and on completions of any integrable representation.

As we showed above, for any given subalgebra $\so_\al$ the KM algebra decomposes 
into a direct sum of $\so_\al$-modules which have Hilbert space completions. More specifically, for each timelike imaginary
root $\al$ (and each choice of a multicommutator $E(\alpha)$ of its root space $\mf{g}_\al$) we have
\beq
 \mf{g}(A)  = (\so_\al \otimes \mathbb{C})\oplus \cH_\al
\eeq
with the associated `total' Hilbert space $\cH_\al$, which itself is the direct sum of infinitely 
many Hilbert spaces corresponding to the irreducible representations of this $\so_\al$. We have written $\mf{g}(A)$ on the left side of the equation above, 
but it should be considered a certain completion. 
The reason why we must now consider the completion of the KM algebra with respect to the norm (\ref{HF}), 
rather than algebraically in the sense of formal sums, is that the operators $J^\pm$ and $J_3$ 
are no longer locally nilpotent. The key issue here is the fact that 
these operators and their exponentials are {\em unbounded operators}.\footnote{As is
  well known the usual operator norm does not exist for unbounded operators. However,
  this is not in contradiction with (\ref{Jnorm}) because the norm of the Lie
  algebra elements induced by the bilinear form (\ref{HF}) is different
  from the   standard operator norm.}
This means that they can be defined only on a domain, 
that is, a {\em dense subspace} of any given Hilbert space. 
Dense subspaces are, for instance, obtained by considering finite linear
combinations of the basis elements of the given representation space.  
But then one faces the problem that repeated action of different exponentials 
on any element of the dense subspace may throw one out of the domain, so
a group action (corresponding to a repeated application of exponentials) is not possible in general.

The key idea which allows us to circumvent the difficulty with unbounded
operators is to exponentiate the (adjoint) action of these operators on 
the KM algebra by exploiting our knowledge of $SO(2,1)$ unitary irreducible representations (UIRs), 
and the fact that this group action is defined on the full Hilbert space,
and not just on dense subspaces.
For the exponentiation of  linear combinations of $J^+$, $J^-$ and $J_3$ 
in $\so_\al$ this requires using the unitary operators
\begin{align}
\label{eq:uexp}
U(w,r) \,=\,
\exp\big( iw J^+ + i\overline{w} J^- + ir  J_3 \big) \quad
\Rightarrow \quad U^\dagger = U^{-1}
\end{align}
with $w\in\cx$ and $r \in \reals$. This action corresponds to the 
compact real form of $\mf{g}(A)$, with anti-\Herm  generators
\begin{align}
\label{eq:cptg}
i(J^+ + J^-) \;, \quad (J^+ - J^-) \;,\quad iJ_3
\end{align}
so this is {\em not} the exponential of the standard split real form. As the unitary representation 
spaces of $SO(2,1)$ are complex, the adjoint action here is on the complexified
Lie algebra $\mf{g}(A)_\cx$.

On each of these representation spaces, hence on all of $\cH_\al$ and the
KM algebra $\mf{g}(A)$, the action of the  {\em group} $SO(2,1)_\al$ or a cover is well defined. 
Concretely this is done as follows: for all UIRs we can evaluate the the action of $SO(2,1)$
on any basis vector $v_n$ by exploiting the action (\ref{Mobiusf}) to express the
transformed function again in terms of the chosen basis:
\beq
S \circ v_m = \sum_{n\in\ints}  {\rm U}_{mn}(S) v_n
\eeq
for any given M\"obius transformation (or in a cover of the M\"obius group). The infinite 
unitary matrix U$_{mn}(S)$ here is the same as in (\ref{S}).  To determine this
matrix for any given unitary operator (\ref{eq:uexp}) 
we use formula (\ref{realexpression}) to convert the argument of the exponential  into a real expression
in terms of the real $\mf{sl}(2)$ generators $\{e,f,h\}$ with real parameters $\{u,v,r\}$.
With these data we can  then compute the two-by-two matrix $S$ and the transformation (\ref{Mobiusf}),
which yields the coefficients U$_{mn} =$ U$_{mn}(S)$  upon expansion in the 
appropriate basis of the relevant function space,
at least in principle. Conversely, given a matrix $S$ we can re-express it in
Iwasawa form and then convert each factor in the Iwasawa decomposition
back to the complex form by  reading (\ref{realexpression}) from right to left.

For the action on the KM algebra we simply replace
the basis vectors by the corresponding elements of the KM algebra. Because the action is on
complex functions the coefficients U$_{mn}(S)$ are in general complex, hence we have
to implement the action of the complexified KM algebra $\mf{g}(A)_\cx$. In this way, we have 
therefore succeeded in exponentiating the associated imaginary root generator.
The presence of non-trivial covers is reflected in the periodicity
properties of (\ref{eq:uexp}) with respect to the rotation parameter $r$: for a $k$-fold covering 
of $SO(2,1)$ we have $r\in \reals/2\pi k \,\nats$.

Requiring the exponential of $\so_\al$ to belong to the KM group therefore requires two generalizations of the common definition of KM group (either minimal or completed). 
The first one is that having both $E(\alpha)$ and $F(\alpha)$ to have well-defined exponentials requires a completion in both Borel directions. The second generalization is 
that we also have to consider covers. The `sheetedness' of the cover is given by $\alpha^2$; as $\alpha^2$ varies we will exhaust all possible covers of $SO(2,1)$.
We therefore find a Kac--Moody group structure that involves the \textit{universal} cover of $SO(2,1)$. 

In physics applications, the physical system is typically built on the symmetric space $G/K$ of a (split real) Kac--Moody group divided by its maximal compact subgroup. Here, one has to make a choice of which Kac--Moody group to use. Using an Iwasawa parametrization of the symmetric space for the minimal group~\cite{DeMedts:2009}, a natural choice uses only a Borel subgroup with Lie algebra corresponding to the positive roots, both real and imaginary, of the KM algebra. The biggest group that can be associated with this Lie algebra is the maximal KM group using democratically one-parameter subgroups associated with \textit{all} generators belonging to positive roots. One parametrization of the symmetric space $G/K$ for the maximal KM group $G$ can then be obtained using for example the standard form parametrization of~\cite[Thm.~8.51]{Marquis} although a different parametrization is common in physics. That this choice of Kac--Moody group is consistent with applications has been discussed for example in~\cite{Bossard:2021jix,Bossard:2021ebg}. Whether our more unitary choice is useful in physics remains to be seen.

\section{\texorpdfstring{Example: The Rank 2 Fibonacci algebra $Fib$}{Example: The Rank 2 Fibonacci algebra Fib}}
\label{sec:Fibadj}

\begin{figure}[h!]
\centering
\includegraphics[scale=1]{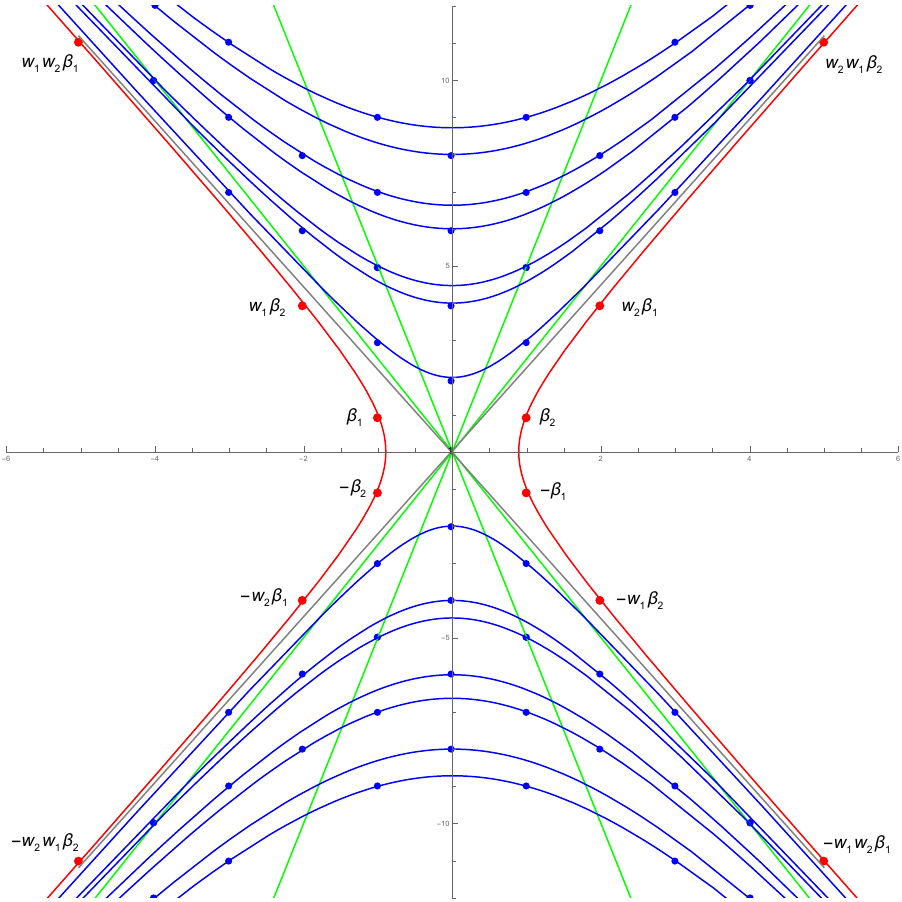}
\caption{\label{fig:Fib}Partial $Fib$ root system with real roots labelled and hyperbolas of constant square length shown. 
First published in \cite{Carbone-Feingold-Freyn:2019}, where the simple roots were labelled by $\beta_i$ instead of $\alpha_i$.}
\end{figure}

In this section we discuss two  kinds of decomposition for the specific example when the Cartan matrix is 
$2\times 2$ with $a_{12} = -3 = a_{21}$, in which case the KM algebra $\mf{g}(A)$ is called $Fib$ 
because its real roots can be described by the Fibonacci numbers~\cite{Feingold:1980}. 
Figure~\ref{fig:Fibpos}, taken from \cite{Feingold-Nicolai:2004}, 
shows a diagram of some of the positive roots for $Fib$, including root multiplicities for the imaginary
roots (black dots). The open circles show some of the positive real roots, with the simple roots,
$\al_1$ and $\al_2$ at the left of the diagram. The central vertical line is the symmetry line of
the outer automorphism which switches $\al_1$ and $\al_2$, but the two angled black lines are the
lines fixed by the two simple Weyl group reflections, $w_1$ and $w_2$. The Weyl group $W$ for $Fib$ is
the infinite dihedral group $D_\infty = \la w_1, w_2\mid |w_1 w_2| = \infty\ra$. 
Figure~\ref{fig:Fib}, taken from \cite{Carbone-Feingold-Freyn:2019}, shows some 
of the roots of $Fib$, both positive and negative, with the real roots on the red hyperbola labelled as 
Weyl conjugates of the simple roots, and the imaginary roots on blue hyperbolas in the light cone. 
The gray lines are the asymptotes of these hyperbolas, the null cone of zero norm points, but no roots
of $Fib$ have zero norm. The inner-most green lines are the fixed lines of the simple Weyl reflections,
and other green lines are their images under the Weyl group action. In each half of the light cone, the 
wedge between the inner green lines is a fundamental domain for the action of $W$ on that half of the 
lightcone, which is tessellated by $W$. 

\subsection{\texorpdfstring{Decomposition of $Fib$ with respect to a real simple root}{Decomposition of Fib with respect to a real simple root}}
\label{sec:Fibreal}

First we will discuss in some detail a (partial) decomposition of $Fib$ into finite-dimensional 
modules with respect to the 
subalgebra $\mf{sl}(2)_1$ with basis $\{e_1, f_1, h_1\}$. The following decomposition can be performed for either the complex Lie algebra or the split real form.
Such a type of decomposition has been used frequently, especially in the physics literature~\cite{Damour:2002cu,Nicolai:2003fw,West:2002jj,Kleinschmidt:2003mf}, 
as well as in the math literature~\cite{Benkart-Kang-Misra:1993}, 
but we present it here for completeness and also for comparison with the case of the imaginary root that is discussed in section~\ref{sec:Fibimag} 
and the one of central interest in this work.

\begin{figure}[h!]
\centering
\includegraphics[angle=90,scale=.70]{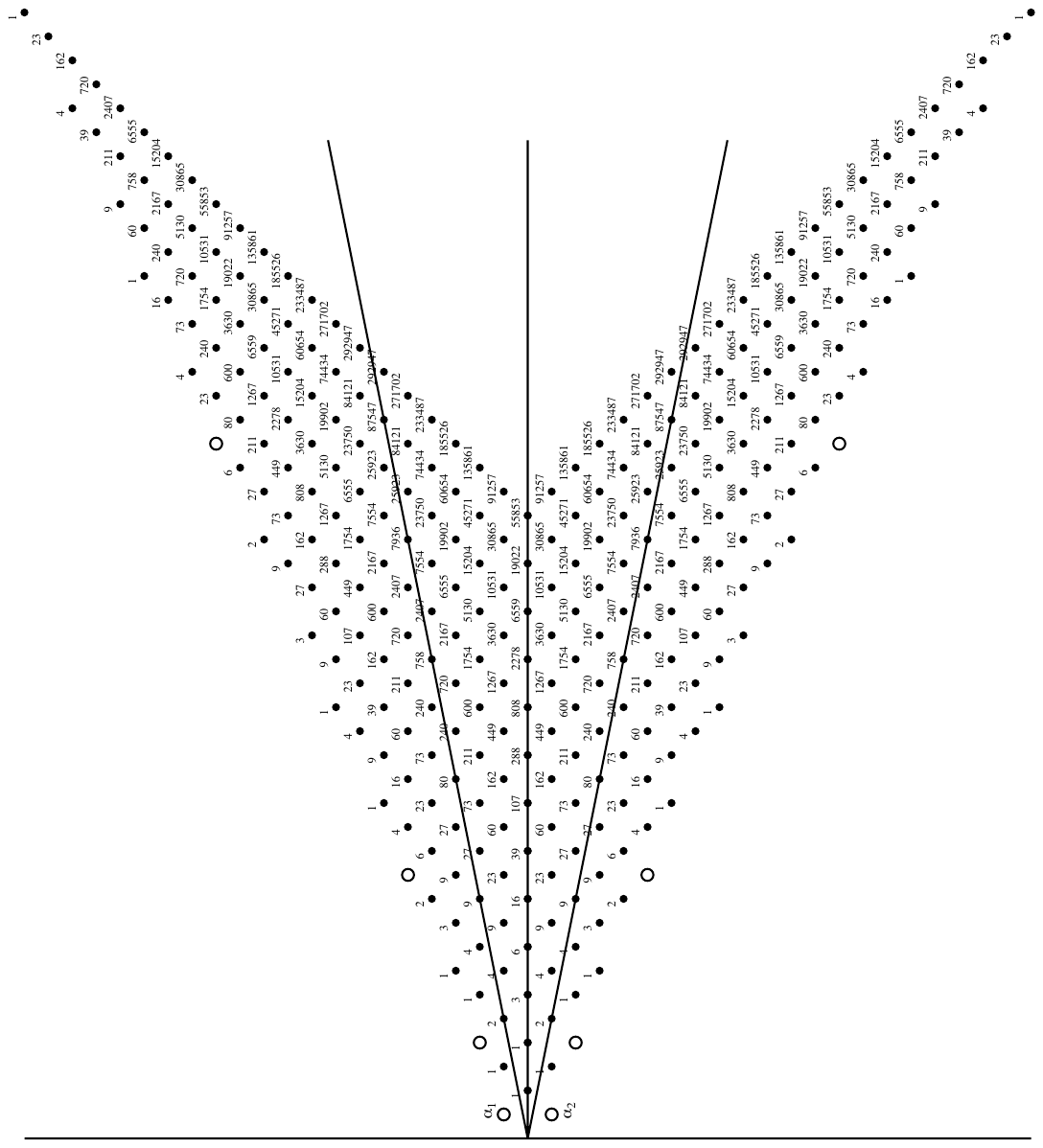}
\caption{Some positive roots of $Fib$ with root multiplicities. First published in \cite{Feingold-Nicolai:2004}.\label{fig:Fibpos}}
\end{figure}

In Figure~\ref{fig:Fibpos}, the subalgebra $\mf{sl}(2)_1$ corresponds to the simple root $\al_1$, giving a direction for the decomposition of $Fib$ into finite-dimensional irreducible 
$\mf{sl}(2)$-modules. For $0\leq m\in\ints$, we denote by $V(m)$ the irreducible $\mf{sl}(2)$-module with $\dim(V(m)) = m+1$. Since $Fib$ has a 
symmetric Cartan matrix, the two choices for simple root, $\al_i$, yield symmetric decompositions, so it is enough to just look at one choice. 
Of the two open circles corresponding to simple roots, let the one towards the left be $\al_1$, and the one to the right be $\al_2$. In the decomposition
of $Fib$ with respect to $\mf{sl}(2)_1$, the first irreducible representation (irrep) to note is $\mf{sl}(2)_1$ itself, a copy of $V(2)$. Since the 
dimension of the Cartan subalgebra (corresponding to the origin in Figure~\ref{fig:Fibpos}) is $2$, there should be another $\mf{sl}(2)_1$ irrep with a weight space
in that Cartan. In fact, it is trivial to compute $[e_1, ah_1+bh_2] = (-2a+3b) e_1$ and $[f_1, ah_1+bh_2] = (2a-3b) f_1$, so that for $2a=3b$ we find
the one-dimensional span of $3h_1+2h_2$ is a trivial module $V(0)$, and no other irrep has a non-trivial intersection with the Cartan subalgebra. 
The next irrep is the one generated by $e_2$ with basis 
\begin{equation}
\{e_2, \ \ [e_1,e_2] = e_{12},\ \  [e_1,[e_1,e_2]] = e_{112},\ \  [e_1,[e_1,[e_1,e_2]]] = e_{1112}\}
\end{equation}
where we have used the notation for multibrackets from section~\ref{sec:KMAs}. This irrep $V(3)$ corresponds to the root string $\al_2$, $\al_2+\al_1$, $\al_2+2\al_1$, 
$\al_2+3\al_1$ in the root diagram, whose end points are real roots and whose middle points are imaginary roots, each with multiplicity $1$. There 
can be no other irreps with weights on that line of roots, so we go to the next parallel line starting with the imaginary root $\al_1+2\al_2$. That root space
is one-dimensional with basis vector $e_{212}$ since there is no way to get to that root space except by $[e_2,[e_1,e_2]]$. The $\mf{sl}(2)_1$ irrep 
generated by that root space is $5$-dimensional because the Weyl group reflection $w_1$ sends $\al_1+2\al_2$ to 
\begin{equation}
w_1(\al_1+2\al_2) = -\al_1 + 2(\al_2+3\al_1) = 5\al_1 + 2\al_2
\end{equation}
forcing the weights of the irrep $V(4)$ to be the $\al_1$-string $\{m\al_1+2\al_2\mid 1\leq m\leq 5\}$ and a basis for that $V(4)$ must be
\begin{equation}
\{ e_{212}, e_{1212}, e_{11212}, e_{111212}, e_{1111212} \}.
\end{equation}
But the dimension of the root space $3\al_1+2\al_2$ is $2$ so there must also be a trivial $1$-dimensional module in that root space. To find an
explicit basis vector for it, take an arbitrary linear combination of the two independent multibrackets in that root space, $e_{11212}$ and $e_{21112}$,
and solve 
\begin{equation}
[f_1, a e_{11212} + b e_{21112}] = 0.
\end{equation}
The same linear condition will occur if $f_1$ is replaced by $e_1$, since this is a trivial module, and it will give a basis for a $V(0)$ irrep in that root space. 
Those two irreps fill up that root string, so the decomposition process continues on the next parallel $\al_1$-string of roots starting with the real root $\al_1+3\al_2$,
$\{m\al_1+3\al_2\mid 1\leq m\leq 8\}$. The list of root multiplicities from Figure~\ref{fig:Fibpos} for that string is $\{1, 2, 3, 4, 4, 3, 2, 1\}$ so we see that there must
be an irrep $V(7)$ having all eight of those weights, as well as an irrep $V(5)$ having the middle six weights, plus an irrep $V(3)$ having the middle
four weights, plus an irrep $V(1)$ having the middle two weights. No other irreps occur in that string, and each one occurred only once. It would be
straight-forward but tedious to find explicit basis vectors for those irreps, or just lowest weight vectors. 

There would be no difference in the above decomposition if we were looking at the split real form $Fib_{\reals}$, but each irrep would be a real vector space.
Some interesting patterns have been seen in such a $\ints$-graded decomposition where the grading of the irreps is according to the coefficient of $\al_2$ in
the $\al_1$-string. The obvious symmetry between positive and negative graded pieces means it suffices to understand the positively graded part. The Lie 
bracket respects the grading, of course, and the $0$-graded piece is just $\mf{sl}(2)_1$ plus the trivial module in the Cartan. So the idea is to see how the
irrep $V(3)$ comprising the $1$-graded piece, bracketed with itself is related to the $2$-graded piece, which is the sum $V(4)\oplus V(0)$. Naturally, the
bracket should correspond to anti-symmetric tensors in the tensor product $V(3)\otimes V(3)$, and we have complete information about such a tensor
product decomposition of finite-dimensional irreps of $\mf{sl}(2)$ from the theory of Clebsch--Gordan. 
One finds that the wedge product $V(3)\wedge V(3)$
exactly equals the sum $V(4)\oplus V(0)$. One expects to get the $3$-graded piece by bracketing the $1$-graded piece $V(3)$ 
with the $2$-graded piece, and so on recursively. But this expansion seems to just get more and more complicated as the grading increases, with no clear
pattern emerging. A similar situation was encountered \cite{Feingold-Frenkel:1983} in the decomposition of the rank-three hyperbolic $\cF$  (mentioned in option (2) in the introduction) 
with respect to its $A_1^{(1)}$ affine subalgebra, where the $\ints$-grading
was with respect to the ``level'' of the affine submodules. A clear answer for level $2$ gave a closed generating function for infinitely many imaginary roots of
$\cF$ because level $1$ was a single irrep whose multiplicities were exactly the values of the classical partition function. Higher levels were studied in \cite{Kang:1994}, up to level 4, and in \cite{Bauer-Bernard:1997} up to level 3, but this method has never yielded a new insight into the full structure of a hyperbolic KM algebra. 

In Figure~\ref{fig:Vrho} we have a graphical display of some weights of the irreducible highest weight $Fib$-module with highest weight $\rho = \lambda_1 + \lambda_2$
along with the weight multiplicities. The weights are determined by 
the action of the Weyl group $W$ 
and the root-string properties of finite-dimensional $\mf{sl}(2)$-modules. The multiplicities are determined recursively by the
Racah--Speiser formula \cite{Kolman-Beck:1973}, which is valid in any irreducible highest weight module $V^\lambda$:
\begin{equation}
\text{Mult}_\lambda(\mu) = \sum_{1\neq w\in W} (-1)^{\ell(w)+1} \text{Mult}_\lambda(\mu + \rho - w(\rho)).
\end{equation}
This is valid for any weight $\mu$ of $V^\lambda$ not in the Weyl orbit of highest weight $\lambda$. All weights in the Weyl orbit $W\cdot\lambda$ have 
multiplicity $1$. For $Fib$ it is easy to compute $\rho - w(\rho)$ and see that for $w\in\{w_1, w_2, w_1w_2, w_2w_1, w_2w_1w_2, w_1w_2w_1\}$ it equals
\begin{equation}
\al_1, \quad \al_2, \quad \al_1+4\al_2, \quad 4\al_1+\al_2, \quad 12\al_1+4\al_2, \quad 4\al_1+12\al_2.
\end{equation}
For the weights of $V^\rho$ shown in Figure~\ref{fig:Vrho}, this recursion only needed Weyl group elements of length $\ell(w)\leq 2$, so only the first four shifts in the
last list. The same algorithm was applied to the fundamental representation $V^{\lambda_1}$ and the results are shown in Figure~\ref{fig:Vlam1}. 
The point of displaying those weight diagrams is to help understand how those modules decompose under the action of  subalgebras of 
$Fib$ like the root subalgebra $\mf{sl}(2)_1$ or an $\so_\al$ subalgebra for an imaginary root $\alpha$ that we discuss in the next section.

\subsection{\texorpdfstring{Decomposition of $Fib$ with respect to an imaginary subalgebra}{Decomposition of Fib with respect to an imaginary subalgebra}}
\label{sec:Fibimag}

Let $e_i$, $f_i$, $h_i$ for $i = 1,2$, be the generators of $Fib$ with the KM relations 
coming from the Cartan matrix $A$. The positive imaginary root of $Fib$ with lowest height is $\al = \al_1+\al_2$,
corresponding the the multi-bracket $e_{12} = - e_{21}$, which we take as $E(\al)$ in the sense of section~\ref{sec:so21s}. Then $F(\al) = -\omega(E(\al)) =
- f_{12} = f_{21}$, and $H(\al) = h_1 + h_2$. We find that $\al^2 = -2$, and using the Jacobi identity, 
we compute the bracket
\begin{align}
[E(\al),F(\al)] &= [e_{12}, f_{21}] = [[e_{12},f_2], f_1] + [f_2, [e_{12},f_1]] = [[e_1,h_2], f_1] + [f_2, [h_1,e_2]]\nn \\
&= [3e_1, f_1] + [f_2, -3e_2] = 3(h_1 + h_2)
\end{align}
so $N = 3$ in formula (\ref{EFH}). This gives us the basis of $ \so_\al$ as in (\ref{JJJ}),
\beq\label{JJJFib}
J_3 = -\frac{h_1 + h_2}{2} \;,\quad 
J^+ = \frac{e_{12}}{\sqrt{6}} \;,\quad J^- = \frac{f_{21}}{\sqrt{6}} \;\;.\quad
\eeq
We begin to find the decomposition of $Fib$ into a direct sum of irreducible $\so_\al$-modules, where the action
of $\so_\al$ is the adjoint action in $Fib$. It is easy to check that 
\begin{equation}
0 = [J_3, h_1 - h_2] =  [J^\pm, h_1 - h_2] 
\end{equation}
so the one-dimensional subspace spanned by $v_0 = h_1 - h_2$ is a trivial $\so_\al$-module. 
Note that for any root $\beta = n_1\al_1 + n_2\al_2$ we have $\beta(h_1 - h_2) = 5(n_1 - n_2)$, so for any 
root vector $x_\beta\in\mf{g}_\beta$, $\frac{1}{5} [h_1-h_2, x_\beta] = (n_1 - n_2) x_\beta$ so the operator 
$\frac{1}{5} (h_1 - h_2)$ provides a $\ints$-grading on $Fib$ corresponding to the horizontal position of the root $\beta$ in Figure~\ref{fig:Fibpos}. 

Since $J_3$ and $v_0$ are independent, they form a basis for the Cartan subalgebra of $Fib$, so there cannot
be any other irreducible $\so_\al$-modules in the decomposition having a non-trivial intersection with the Cartan. 
In particular, this means that only discrete series modules (highest or lowest weight modules) can occur on the
central symmetry line of roots $n\al = n(\al_1+\al_2)$. Looking at parallel lines just to the side shifted by either
adding $\al_1$ or $\al_2$, let us see what $\so_\al$-modules are generated by the simple root vectors $e_1$ and $e_2$.
We find that for $i = 1,2$, 
\begin{align}
[J_3, e_i] &= -\frac{1}{2} [h_1 + h_2, e_i] = -\frac{1}{2} \al_i(h_1 + h_2) e_i  = \frac{1}{2} e_i,\nn\\
[J^+, e_1] &= -\frac{1}{\sqrt{6}} e_{112}, \quad [J^-, e_1] = \frac{3}{\sqrt{6}} f_{2}, \quad 
[J^+, e_2] = -\frac{1}{\sqrt{6}} e_{212}, \quad [J^-, e_2] = -\frac{3}{\sqrt{6}} f_{1}
\end{align}
and then, using the Jacobi identity, we get 
\begin{equation}
[J^-,[J^+,e_1]] =  -\frac{1}{6} [f_{21}, e_{112}] = 2 e_1\quad\text{and}\quad 
[J^+,[J^-,e_1]] =  \frac{3}{2} e_1
\end{equation}
so the Casimir operator on $e_1$ gives
\begin{equation}
\Omega e_1 = (J_3 J_3 - J^- J^+ - J^+ J^-) e_1 
= -\frac{13}{4} e_1.
\end{equation}
A similar calculation for $e_2$ (or using the symmetry exchanging subscripts $1$ and $2$), gives 
\begin{equation}
[J^-,[J^+,e_2]] =  -\frac{1}{6} [f_{21}, e_{212}] = 2 e_2\quad\hbox{and}\quad 
[J^+,[J^-,e_2]] =  \frac{3}{2} e_2
\end{equation}
so 
\begin{equation}
\Omega e_2 = -\frac{13}{4} e_2.
\end{equation}
This means there are two principal series $\so_\al$-modules generated by these two simple root vectors, one with weights
in the line of roots $\al_1 + n\al$, and the other in the line of roots $\al_2 + n\al$. In both cases we find the parameter $s$ such that $s(s-1) = -\frac{13}{4}$ to be
\begin{equation}
s = \frac{1\pm i q}{2} 
\quad \hbox{with}\quad q = \sqrt{12}.
\end{equation}

Going back to the center line of symmetry, Figure~\ref{fig:Fibpos} shows that for $\beta = 2(\al_1 + \al_2)$, the $\beta$ root space has
dimension $1$, and it is easy to check that $e_{1212} = e_{2112}$ is a basis for it. To verify that it is a lowest weight vector
killed by $J^-$ we use that 
\begin{equation}
[f_1, e_{1212}] = 4 e_{212}\quad \hbox{and}\quad [f_2, e_{1212}] = 4 e_{112}
\end{equation}
as well as 
\begin{equation}
[f_1, e_{112}] = 4 e_{12} \quad\hbox{and}\quad [f_2, e_{212}] = 4 e_{12} 
\end{equation}
to compute
\begin{equation}
[f_{12}, e_{1212}] =  [[f_1, e_{1212}], f_2] + [f_1, [f_2,e_{1212}]] = 4 [e_{212}, f_2] + 4 [f_1, e_{112}] = 0.
\end{equation}
Since $\beta(h_1+h_2) = -4$ we have
\begin{equation}
[J_3, e_{1212}] = -\frac{1}{2} [h_1+h_2, e_{1212}] = -\frac{1}{2} \beta(h_1+h_2) e_{1212} = 2 e_{1212}
\end{equation}
so the parameter $s = 2$ and 
\begin{equation}
\Omega e_{1212} = s(s-1) e_{1212} = 2 e_{1212}.
\end{equation}

Staying on the center line the next space reached by $J^+$ when acting on $e_{1212}$ is the root space of $\beta=3(\alpha_1+\alpha_2)$ which has dimension $3$ and so there should be two lowest weight vectors in that space. A basis of the $3(\alpha_1+\alpha_2)$ root space is given by
\begin{align}
e_{112212}\,,\quad e_{121212}\quad \text{and}\quad e_{211212}\,.
\end{align}
Acting with $J^+$ on $e_{1212}$ leads to 
\begin{align}
J^+ e_{1212} = \frac{1}{\sqrt{6}} \left( e_{121212} - e_{211212} \right)\,.
\end{align}
One can check that the following are independent lowest weight vectors for the action of $J^-$
\begin{align}
\ell_1 = e_{121212}+e_{211212} \quad\text{and}\quad  \ell_2 = e_{112212} + 3e_{211212}\,.
\end{align}
The corresponding Casimir eigenvalues are given from their $J_3$ eigenvalues:
\begin{align}
J_3 \begin{pmatrix}\ell_1\\\ell_2 \end{pmatrix} = \begin{pmatrix} 3 & 0 \\0 & 3\end{pmatrix} \begin{pmatrix}\ell_1\\\ell_2 \end{pmatrix} \,,
\end{align}
so that in both cases $\Omega= 6$.

Let us also determine a lowest weight representation off the center line. There must be one in the root space of $3\alpha_1+2\alpha_2$ that has the basis
\begin{align}
e_{11212} \quad \text{and}\quad e_{21112}\,.
\end{align}
The lowest weight combination annihilated by $J^-$ is
\begin{align}
3 e_{11212} + 4e_{21112}
\end{align}
with $J_3$ eigenvalue $s=\tfrac52$ and so $\Omega= \tfrac{15}4$.

The fractional $s$-value for a discrete series also appears in
\begin{align}
\exp( i r J_3) e_{11212} = \exp( 5 i r/ 2) e_{11212}\,,
\end{align}
so that $r\in \mathbb{R}/ (4\pi \mathbb{Z})$, showing that we are dealing with a cover of $SO(2,1)$.

\subsection{\texorpdfstring{Decomposition of highest weight $Fib$-representations}{Decomposition of highest weight Fib-representations}}

Up until now, we have discussed the decomposition of the KM algebra, $\mf{g}_\cx(A)$, with respect to an imaginary subalgebra, $\so_\al$. 
But one also has the decomposition of any highest or lowest weight representation of the KM algebra with respect to the action
of that imaginary subalgebra. In that situation only discrete series of $\so_\al$ representations can occur in the decomposition, so some
of the complications coming from the continuous series do not arise. 

\begin{figure}[h!]
\centering
\includegraphics[width=5.8in,height=7.5in]{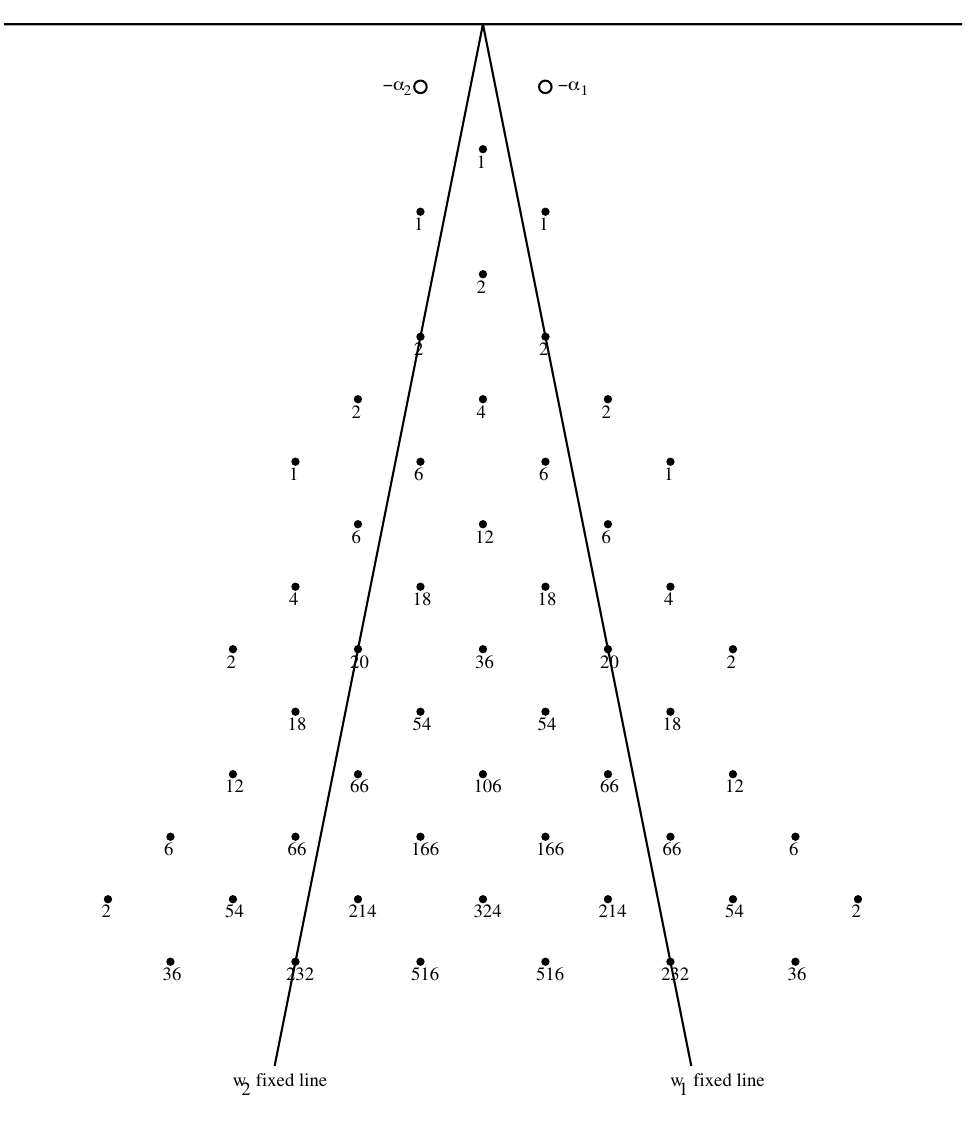}
\caption{\label{fig:Vrho}Some weights of irreducible $Fib$ module $V^\rho$ with highest weight $\rho = \lambda_1 + \lambda_2$ showing multiplicities.}
\end{figure}

\begin{figure}[h!]
\centering
\includegraphics[width=5.8in,height=7.5in]{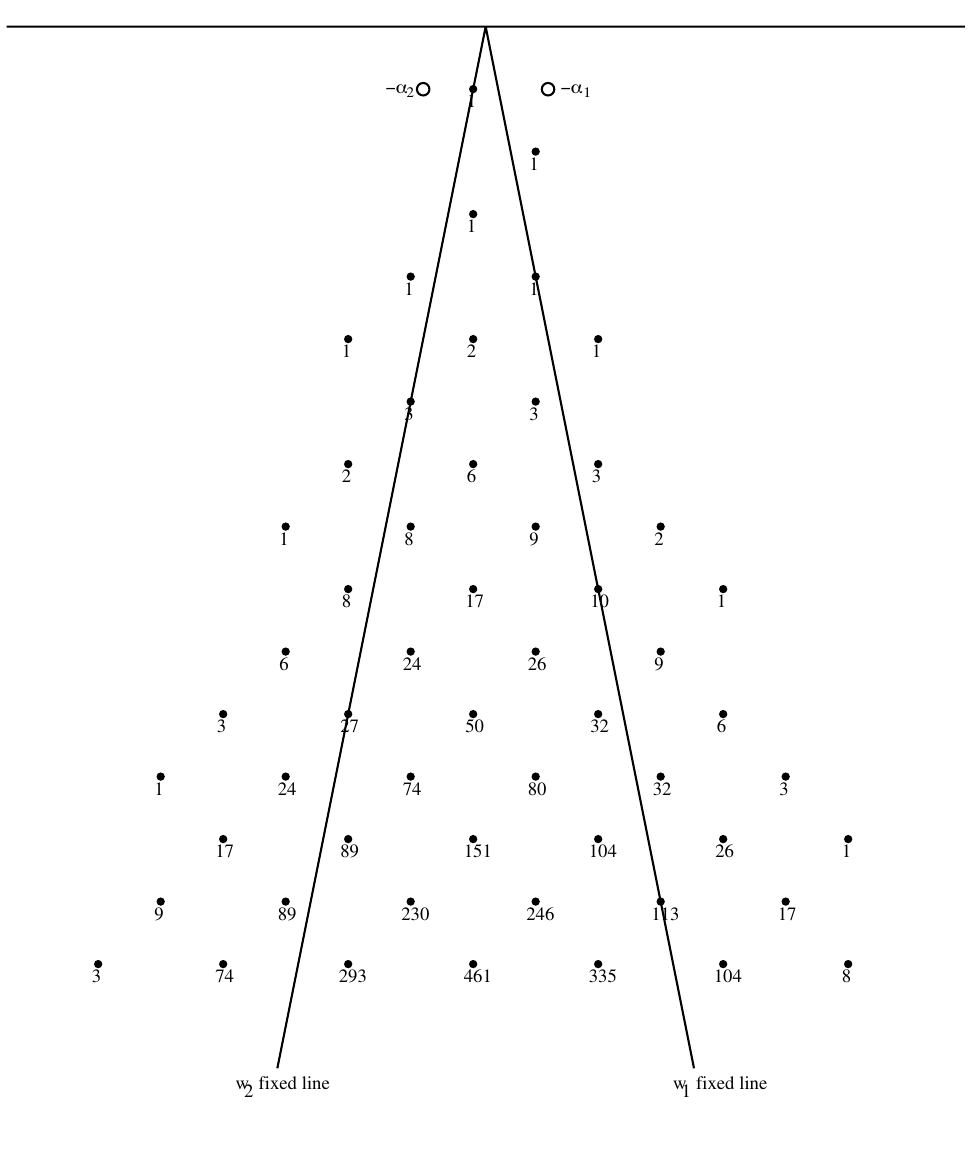}
\caption{\label{fig:Vlam1}Some weights of irreducible $Fib$ fundamental module $V^{\lambda_1}$ with highest weight $\lambda_1$ showing multiplicities.}
\end{figure}

As an illustration, we present here two examples of highest weight 
representations of the rank $2$ hyperbolic KM algebra, $Fib$, $V^\rho$ and $V^{\lambda_1}$, whose partial weight diagrams are shown in
Figure~\ref{fig:Vrho} and Figure~\ref{fig:Vlam1}. In each case, let $v_\lambda$ be a highest weight vector of weight $\lambda$ in $V^\lambda$, so that $J^+(v_\lambda) = 0$
and from (\ref{JJJ1}) we get
\begin{equation}
J_3\cdot v_\lambda = \frac{\lambda(H(\al))}{\al\cdot\al} \ v_\lambda = \frac{\lambda\cdot\al}{\al\cdot\al} \ v_\lambda.
\end{equation}
From (\ref{Casimir1}) we have 
\begin{equation}
\Omega (v_\lambda) = (J_3(J_3 + 1) - 2 J^- J^+ ) (v_\lambda) = J_3(J_3 + 1) (v_\lambda) = 
\left(\frac{\lambda\cdot\al}{\al\cdot\al}\right) \left(\frac{\lambda\cdot\al}{\al\cdot\al} + 1 \right) v_\lambda.
\end{equation}
As we did in the previous section, use the positive imaginary root $\al = \al_1+\al_2$ of $Fib$ corresponding the the multi-bracket $E(\al) = e_{12}$, 
giving the formulas in (\ref{JJJFib}). Then we have 
\begin{equation}
J_3(v_\lambda) = \frac{-1}{2} (h_1+h_2)\cdot v_\lambda = \frac{-1}{2} \lambda(h_1+h_2) v_\lambda
\end{equation}
and in particular, 
\begin{equation}
J_3(v_{\lambda_1}) = \frac{-1}{2} \lambda_1(h_1+h_2) v_{\lambda_1} = \frac{-1}{2} v_{\lambda_1}\quad\hbox{ so }\quad s =  \frac{1}{2}
\end{equation}
and 
\begin{equation}
J_3(v_{\rho}) = \frac{-1}{2} \rho(h_1+h_2) v_{\rho} = - v_{\rho}\quad\hbox{ so }\quad s = 1.
\end{equation}
Therefore, 
\begin{equation}
\Omega(v_{\lambda_1})  
= \frac{-1}{4} v_{\lambda_1}
\quad\quad \text{and}\quad\quad
\Omega(v_{\rho}) = 
0.
\end{equation}
The fractional value of $s$ for the lowest weight vector $v_{\lambda_1}$ means that the relevant group acting on this irrep will be a cover of $SO(2,1)$.

Examining the (partial) weight diagrams of these two modules in Figures~\ref{fig:Vrho} and~\ref{fig:Vlam1}, we see that the vertical line of weights going down from the
highest weight contains the discrete series (\ref{DiscReps}) module $\cD_s^-$ for $s =  \frac{1}{2}$ in $V^{\lambda_1}$ and for $s = 1$ in $V^\rho$. 
In that vertical line of weights $\{\lambda_1 - n\al \mid 0\leq n\in\ints\}$ for $V^{\lambda_1}$ the weight multiplicities shown in Figure~\ref{fig:Vlam1} are 
$\{1, 1, 2, 6, 17, 50, 151, 461\}$ corresponding to $0\leq n\leq 7$. Since the weight spaces in $\cD_s^-$ are each $1$-dimensional, the ``top'' 
summand in the decomposition for that line accounts for the first two $1$'s on that list, and decreases each of the following numbers by $1$. So the
next summand is determined by a highest weight vector (killed by $J^+$) of weight $\lambda = \lambda_1 - 2\al$. We do not explicitly compute that 
highest weight vector here, but it is straightforward to find it as a linear combination of basis vectors in that $2$-dimensional weight space of $V^{\lambda_1}$.
Since $\al(h_1+h_2) = -2$, we see that the next summand is a discrete series module with $s =  \frac{5}{2}$, which accounts for one of the dimensions in
each of the list of multiplicities, reducing the list to $\{0, 0, 0, 4, 15, 48, 149, 459\}$. In general, if there are any highest weight vectors in that column with
weight $\lambda_1 - n\al$, the eigenvalue of $J_3$ on such vectors will be 
\begin{equation}
\frac{-1}{2} (\lambda_1 - n\al)(h_1+h_2) 
= - \frac{2n+1}{2}\quad\hbox{so}\quad s = \frac{2n+1}{2}
\end{equation}
is the corresponding value of parameter $s$ for each copy of the discrete series module $\cD_s^-$ at that weight in the decomposition. 
Clearly this decomposition process continues, giving a $4$-dimensional space of 
highest weight vectors with weight $\lambda_1 - 3\al$, and thus, four copies of $\cD_s^-$ with $s =  \frac{7}{2}$, reducing each of the remaining numbers 
by $4$, leaving the list $\{0, 0, 0, 0, 11, 44, 145, 455\}$. There will be $11$ copies of $\cD_s^-$ with $s =  \frac{9}{2}$, and $33$ copies with $s =  \frac{11}{2}$,
and $101$ copies with $s =  \frac{13}{2}$, and $310$ copies with $s =  \frac{15}{2}$, etc. Each column of weights in the diagram has a top weight, and 
each weight below has a multiplicity, so the process above produces a list of summands consisting of copies of $\cD_s^-$ for values of $s$ determined
by the weight. For example, the column to the right of $\lambda_1$ starts with $\lambda_1 - \al_1$ and consists of weights of the form 
$\{\lambda_1 - \al_1 - n\al \mid 0\leq n\in\ints\}$. Since $\al_1(h_1+h_2) = -1$, the eigenvalue of $J_3$ on such weight vectors is 
$\frac{-1}{2} (\lambda_1 - \al_1 - n\al)(h_1+h_2) 
= -(n+1)$ corresponding to $s = n+1$. In Figure~\ref{fig:Vlam1} we see the list of multiplicities
is $\{1, 1, 3, 9, 26, 80, 246\}$, giving a list of multiplicities of discrete series modules for that column as differences. The complete decomposition involves
doing that process for every column in the weight diagram. The reader is invited to carry out part of this process for $V^\rho$ using the multiplicities shown
in Figure~\ref{fig:Vrho}. 

Now we can apply formulas for the action of the group $SO(2,1)_\al$ and its covers on each of the discrete series summands, $\cD_s^-$, in the decomposition of any
highest weight representation $V^{\lambda}$ of $\mf{g}_\cx(A)$. These can be understood as an exponentiation of the imaginary Lie subalgebra, $\so_\al$,
as operators on $V^{\lambda}$. The infinite sums involved can be understood as converging with respect to a \Herm form on $V^{\lambda}$ which has been
defined in \cite{Kac}, uniquely determined by $(v_\lambda, v_\lambda) = 1$ and $(X(v), w) = (v, X^\dagger(w))$ for every $v,w\in V^{\lambda}$ and 
every $X\in \mf{g}_\cx(A)$. Since a highest weight vector, $v_\lambda \in V^{\lambda}$ is only determined up to a scalar, the same is true of the form. 
 
\subsection{\texorpdfstring{Decomposition of $Fib$ with respect to another imaginary subalgebra}{Decomposition of Fib with respect to another imaginary subalgebra}}
\label{sec:Fibimag2}

Another positive imaginary root of $Fib$, not in the Weyl group orbit of $\al_1+\al_2$, is $\al = 2\al_1+3\al_2$, whose root space is $2$-dimensional with basis
$\{e_{21212}, e_{12221} \}$, so we could take either one of these as $E(\al)$ in the sense of section~\ref{sec:so21s}. For this section we choose $E(\al) = e_{21212}$ so
that $F(\al) = -\omega(E(\al)) = f_{21212}$, and $H(\al) = 2h_1 + 3h_2$. We find that $\al^2 = -10$, and using the Jacobi identity, we compute the bracket
\begin{equation}
[E(\al),F(\al)] = [e_{21212}, f_{21212}] = 288 (2h_1 + 3h_2)
\end{equation}
so $N = 288$ in formula (\ref{EFH}). This gives us the basis of $ \so_\al$ 
\beq
J_3 = -\frac{2h_1 + 3h_2}{10} \;,\quad 
J^+ = \frac{e_{21212}}{24\sqrt{5}} \;,\quad J^- = \frac{f_{21212}}{24\sqrt{5}} \;\;.\quad
\eeq
The bracket calculation above used brackets from section \ref{sec:Fibimag} as well as the following bracket:
\begin{equation}
[e_{1212}, f_{1212}] = -96(h_1 + h_2).
\end{equation}

\begin{figure}[h!]
\centering
\includegraphics[angle=90,scale=.70]{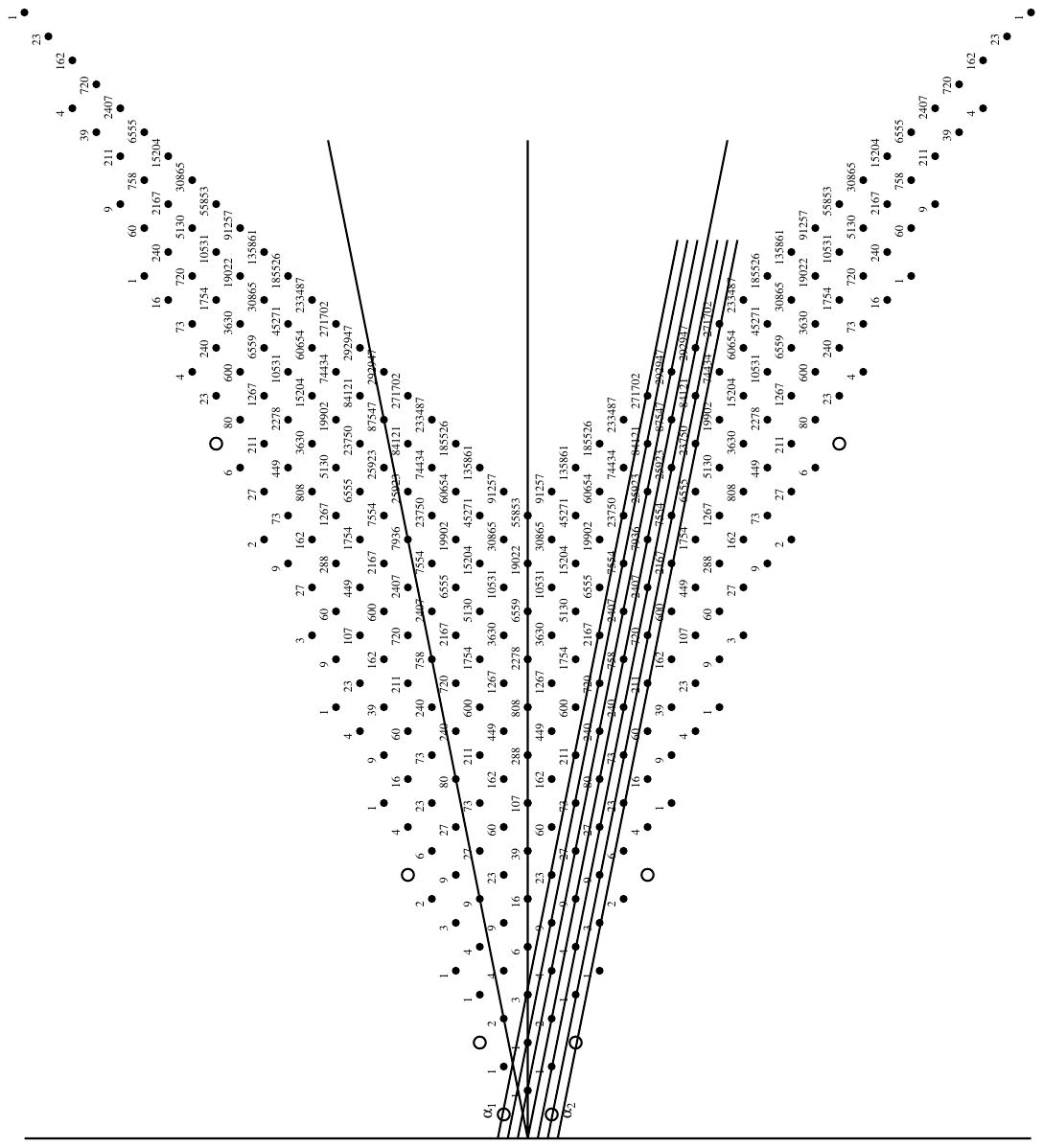}
\caption{Some positive roots of $Fib$ with root multiplicities and six lines parallel to the line through $\alpha = 2\alpha_1+3\alpha_2$ where principal series might occur for
$\so_\al$.\label{fig:Fibpos2}}
\end{figure}

We begin to find the decomposition of $Fib$ into a direct sum of irreducible $\so_\al$-modules whose weights will be on lines parallel to the line through $\al$. 
Since $(2\al_1+3\al_2)(h_2) = 0$, It is easy to see that 
\begin{equation}
0 = [J_3, h_2] =  [J^\pm, h_2] 
\end{equation}
so the one-dimensional subspace spanned by $v_0 = h_2$ is a trivial $\so_\al$-module. 
Since $J_3$ and $v_0$ are independent, they form a basis for the Cartan subalgebra of $Fib$, so there cannot
be any other irreducible $\so_\al$-modules in the decomposition having a non-trivial intersection with the Cartan. 
In particular, this means that only discrete series modules (highest or lowest weight modules) can occur on the
line of roots $n\al = n(2\al_1+3\al_2)$. Looking at Figure \ref{fig:Fibpos2} we see six parallel lines of roots, three on each side of that line, where unbroken $\al$ 
root strings could contain principal series representations. Those lines are each of the form $\{\mu + n\al\mid n\in\ints\}$ for $\mu$ in the set 
\begin{equation}
\{-\al_1, \al_2, \al_1+2\al_2, \al_1+\al_2, -\al_2, \al_1\}
\end{equation}
Let us see what principal series $\so_\al$-modules are generated by root vectors in $\mf{g}_\mu$ for the roots $\mu$ in that set.

For $\mu = -\al_1$ the basis root vector is $f_1$ and for $\mu = -\al_2$ the basis root vector is $f_2$ and we find
\begin{align}
[J_3, f_1] &= -\frac{1}{10} [2h_1 + 3h_2, f_1] = \frac{1}{10} \al_1(2h_1 + 3h_2) f_1  = \frac{-1}{2} f_1,\nn\\
[J_3, f_2] &= -\frac{1}{10} [2h_1 + 3h_2, f_2] = \frac{1}{10} \al_2(2h_1 + 3h_2) f_1  = 0,\nn\\
[J^+, f_1] &= \frac{-1}{6\sqrt{5}} e_{2212}, \qquad [J^-, f_1] = \frac{-1}{24\sqrt{5}} f_{121212},\nn\\
[J^+, f_2] &= \frac{-1}{4\sqrt{5}} e_{1212}, \qquad [J^-, f_2] = \frac{-1}{24\sqrt{5}} f_{221212}
\end{align}
and then, using the Jacobi identity, we get 
\begin{align}
[J^-,[J^+, f_1]] &=  \frac{-1}{720} [f_{21212}, e_{2212}] = \frac{1}{5} f_1,\\
[J^+,[J^-, f_1]] &=  \frac{-1}{2880} [e_{21212}, f_{121212}] = \frac{7}{10} f_1 ,\\
[J^-,[J^+, f_2]] &=  \frac{-1}{480} [f_{21212}, e_{1212}] = \frac{3}{5} f_2,\\
[J^+,[J^-, f_2]] &=  \frac{-1}{2880} [e_{21212}, f_{221212}] = \frac{3}{5} f_2 .
\end{align}
so the Casimir operator on $f_1$ and $f_2$ gives
\begin{align}
\Omega f_1 = (J_3 J_3 - J^- J^+ - J^+ J^-) f_1 &= \left(\frac{1}{4} - \frac{1}{5} - \frac{7}{10}\right) f_1= \frac{-13}{20} f_1 ,\\
\Omega f_2 = (J_3 J_3 - J^- J^+ - J^+ J^-) f_2 &=  \left(0 - \frac{3}{5} - \frac{3}{5}\right) f_2= \frac{-6}{5} f_2 .
\end{align}
Further calculations give the following values of the Casimir operator on basis vectors in $\mu$ root spaces for the other values of $\mu$. 
For $\mu = \al_1$ the basis root vector is $e_1$, for $\mu = \al_2$ the basis root vector is $e_2$, for $\mu = \al_1 + \al_2$ the basis root vector is $e_{12}$
and for $\mu = \al_1 + 2\al_2$ the basis root vector is $e_{212}$. The results are: 
\begin{equation}
\Omega e_1 = \frac{-13}{20} e_1, \qquad
\Omega e_2 = \frac{-6}{5} e_2, \qquad
\Omega e_{12} = \frac{-13}{20} e_{12},\qquad
\Omega e_{212} = \frac{-13}{20} e_{212}.
\end{equation}
Since the Casimir operator values above are all less than $-1/4$, these principal series representations are not complementary, consistent with our
conjecture. Also, note that the Lie algebra automorphism, $\omega$, commutes with the Casimir operator, $\Omega$, so if $\Omega x = \lambda x$
then $\Omega (\omega(x)) = \lambda \omega(x)$. Since $\omega(e_i) = -f_i$, for $i = 1,2$, this explains why the $\Omega$ eigenvalues of $e_i$ and 
$f_i$ are equal. Furthermore, the root $-\al_1 - \al_2$ is in the same $2\al_1 + 3\al_2$ root string as $\al_1 + 2\al_2$, and $\omega(e_{12}) = f_{12}$
is in the $-\al_1 - \al_2$ root space, so that explains why the $\Omega$ eigenvalues of $e_{12}$ and $e_{212}$ are equal.

In the $Fib$ example above and in section~\ref{sec:Fibimag} no complementary series representations arose among the finitely many continuous series $\so_\al$-representations appearing in the decomposition of $Fib$. We have performed similar checks for the rank-two hyperbolic algebra with Cartan matrix $\begin{psmallmatrix}2&-4\\-4&2\end{psmallmatrix}$ as well as the rank-three hyperbolic algebra $\mathcal{F}$ studied in~\cite{Feingold-Frenkel:1983} and found no complementary series representations. Based on this data we make the following conjecture whose further analysis we leave to future work.
\vspace{2mm}

\noindent{\bf Conjecture}: In the decomposition of any hyperbolic KM Lie algebra $\mf{g}$ with respect to $\so_\al$ for positive imaginary root, $\al$, the principal series 
representations that occur are not complementary series.

\newpage

\appendix

\section{\texorpdfstring{Real forms of $\mf{sl}(2,\cx)$}{Real forms of sl(2,C)}}
\label{sec:sl2}

We here summarize some very basic facts about two kinds of real forms of the complex Lie algebra $\mf{sl}(2,\cx)$
of $2\times 2$ complex matrices with trace $0$. This simple $3$-dimensional Lie algebra has basis 
\beq\label{efhbasis}
\left\{e = \bm 0&1\\ 0&0\ebm, f = \bm 0&0\\ 1&0\ebm , h = \bm 1&0\\ 0&-1\ebm  \right\}
\eeq
with the Lie brackets $[e,f] = h$, $[h,e] = 2e$ and $[h,f] = -2f$. Its finite-dimensional representations 
play a crucial role in the representation theory of semi-simple Lie algebras over $\cx$, as well as in the definition of Kac--Moody
Lie algebras. A real form of a complex Lie algebra $\mf{g}_\cx$ is a Lie algebra $\mf{g}$ over $\reals$ such
that $\mf{g}\otimes\cx$ is isomorphic to $\mf{g}_\cx$. For example, the {\it split} real form $\mf{sl}(2,\reals)$ of $\mf{sl}(2,\cx)$ is just
the real span of the basis $\{e,f,h\}$, and it can be understood as the fixed points in $\mf{sl}(2,\cx)$ of the involution 
$\tau$ defined by $\tau(ae) = \ba e$, $\tau(af) = \ba f$ and $\tau(ah) = \ba h$ for any $a\in\cx$. 

An important point is that for any $A\in GL(2,\cx)$, and any real subalgebra $\mf{a}$ of $\mf{sl}(2,\cx)$, the conjugate
$A \mf{a} A^{-1}$ is clearly a real subalgebra of $\mf{sl}(2,\cx)$ isomorphic to $\mf{a}$. For example, if $\mf{a} = \mf{sl}(2,\reals)$,
then each $A \mf{a} A^{-1}$ is a split real form of $\mf{sl}(2,\cx)$, but the entries of its matrices can be complex. It is therefore 
somewhat misleading to speak of ``the split real form'' of $\mf{sl}(2,\cx)$ unless one understands this equivalence of conjugates. 
The same consideration applies to ``the compact real form'' of $\mf{sl}(2,\cx)$, and means that there can be infinitely many 
conjugate versions of it. Below we will give some explicit realizations of these two kinds of real forms. Of course, if $A$ is a real matrix,
the conjugation $A \mf{a} A^{-1}$ amounts to a real change of basis within $\mf{a}$ which means $A \mf{a} A^{-1} = \mf{a}$. 

Two real forms of $\mf{sl}(2,\cx)$ can be found as real Lie algebras of $3\times 3$ matrices. By definition, 
\begin{equation}
\mf{so}(3) = \{A \in\reals_{3\times 3} \mid A^t = -A\}
\end{equation}
is the real Lie algebra of anti-symmetric matrices which has a basis
\beq\label{so3basis}
\left\{M_1 = \bm 0&0&0\\ 0&0&1\\ 0&-1&0 \ebm, M_2 = \bm 0&0&1\\ 0&0&0\\ -1&0&0 \ebm , M_3 = \bm 0&1&0\\ -1&0&0\\ 0&0&0 \ebm  \right\}
\eeq
with the Lie brackets $[M_1,M_2] = M_3$, $[M_2,M_3] = M_1$ and $[M_3,M_1] = M_2$. This is clearly isomorphic to the real Lie algebra $\reals^3$
where the Lie bracket is just the cross product. It is also the real Lie algebra of matrices determined by the standard dot product in $\reals^3$, 
$X\cdot Y = X^t Y$. The adjoint of $A \in\reals_{3\times 3}$ with respect to this dot product is the unique $A^*$ such that $(AX)\cdot Y = X\cdot (A^*Y)$ 
so $(AX)^t Y = X^t A^t Y = X\cdot (A^t Y)$ gives $A^* = A^t = -A$ for $A\in\mf{so}(3)$. 
For any real symmetric $3\times 3$ matrix, $\eta$, we have a symmetric bilinear form
$B_\eta(X,Y) = X^t \eta Y$ and the associated Lie algebra is 
\begin{equation}
\mf{so}(B_\eta) = \{A \in\reals_{3\times 3} \mid A^t\eta = -\eta A\}.
\end{equation}
With respect to $B_\eta(X,Y)$ the adjoint of $A$ is determined by
$B_\eta(AX,Y) = B_\eta(X,A^* Y)$, that is, $(AX)^t \eta Y = X^t \eta (A^* Y)$ so $X^t A^t \eta Y = X^t \eta (A^* Y)$ so $A^*$ is determined by the 
condition $A^t \eta = \eta A^*$. For $A\in \mf{so}(B_\eta)$ this says $A^* = -A$. 
For $\eta = I_3$ this is just $\mf{so}(3)$, but if we use $\eta = diag(-1,1,1)$ we get 
\begin{equation}
\so = \{A \in\reals_{3\times 3} \mid A^t\eta = -\eta A\}
\end{equation}
which has a basis
\beq{\label{so21basis}}
\left\{J_0 = \bm 0&0&0\\ 0&0&1\\ 0&-1&0 \ebm, J_1 = \bm 0&1&0\\ 1&0&0\\ 0&0&0 \ebm, J_2 = \bm 0&0&1\\ 0&0&0\\ 1&0&0 \ebm \right\}
\eeq
with the Lie brackets $[J_0,J_1] = -J_2$, $[J_1,J_2] = J_0$ and $[J_2,J_0] = -J_1$. We wish to show that the real Lie algebras $\so$ and
$\mf{sl}(2,\reals)$ are isomorphic. First note that the following real linear combinations  
\begin{equation}
J_3 = -J_2 \qquad \text{and}\qquad J^\pm = \frac{1}{\sqrt{2}} (J_0 \pm J_1)
\end{equation}
satisfy the bracket relations
\begin{equation}
[J_3, J^\pm] = \pm J^\pm \qquad \text{and}\qquad [J^+, J^-] = - J_3.
\end{equation}
Then the following elements satisfy the bracket relations for the standard basis of $\mf{sl}(2,\reals)$:
\begin{equation}
H = 2J_3 = -2 J_2,\qquad E = \sqrt{2} J^+ = J_0 + J_1  \qquad \text{and}\qquad F = - \sqrt{2} J^- = -J_0 + J_1.
\end{equation}

Going back to the basis (\ref{so3basis}) of $\mf{so}(3)$, the following complex linear combinations
\begin{equation}
M^\pm = \bi (M_1 \pm \bi M_2)  \qquad \text{and}\qquad M_z = \bi M_3
\end{equation}
satisfy the bracket relations:
\begin{equation}
[M_z, M^\pm] = \pm M^\pm  \qquad \text{and}\qquad [M^+, M^-] = 2 M_z.
\end{equation}
The slight rescaling 
\begin{equation}
E = M^+, \qquad F = M^-  \qquad \text{and}\qquad H = 2M_z
\end{equation}
gives the bracket relations for the standard basis of $\mf{sl}(2,\reals)$. The {\it compact} real form $\mf{so}(3)$ and the split real form 
$\mf{sl}(2,\reals)$ are not isomorphic as real Lie algebras, but they are related by the complex linear map above. 
A simple real transformation of the basis $\{M^+, M^-, M_z \}$ gives the bracket relations of $\so$. 
We have the adjoints $M_i^* = -M_i$ for $1\leq i\leq 3$, and since the standard complex dot product is sesquilinear, complex coefficients
get conjugated in an adjoint, so we get
\begin{equation}
M_z^* = M_z  \qquad \text{and}\qquad (M^\pm)^* = M^\mp.
\end{equation}

\section{Abstract algebra: unitary representations}
\label{sec:absrep}

This appendix contains a brief review of unitary irreducible representations of the groups $SO(2,1)$ and $SL(2,\reals)$ (and further covers) for the readers convenience. Some general references are~\cite{Bargmann:1946me,Pukanszky:1964,Sally:1965,Lang:1985,Howe:1992}, see also the recent thesis~\cite{Hauser:2023} that discusses automorphic aspects.
In appendix~\ref{app:FR}, we present functional realizations of these abstract representations.

The starting point is the Lie algebra $\so$ with commutation relations given by (\ref{SO2,1}). We are looking for unitary irreducible representations, that is complex vector spaces admitting an inner product such that the Lie algebra generators satisfy the hermiticity properties in (\ref{Herm}).
The $\so$ Casimir operator, which commutes with $\so$, is given by
\begin{align}\label{Casimir1}
\Omega\,  &=\,  J_3 J_3 - J^+ J^- - J^- J^+  \nonumber\\[1mm]
              &= \, J_3(J_3 -1) - 2 J^+ J^-  \,=\, J_3(J_3 +1) -   2 J^- J^+\,.
\end{align}
Since the representation is irreducible, by Schur's lemma, $\Omega$ is a constant scalar operator on the 
entire irreducible representation. From these relations we immediately obtain
\begin{align}\label{Casimir2}
2 J^+ J^-   & = J_3(J_3 -1) - \Omega \;\;,\quad
2 J^- J^+    = J_3(J_3 +1) - \Omega  
\end{align}
for any irreducible representation. 

\subsection{Discrete series}

Let us first consider the discrete series representations that admit a highest or lowest weight. Since the discussions are fully analogous in the two cases, we restrict mainly to highest weight representations.

Given a highest weight state $w_0$ we have $J^+ w_0 =0$ and we denote its $J_3$-eigenvalue by
$J_3\, w_0 = - s w_0$. We assume that there is an Hermitian inner product $(\cdot, \cdot)$ on the representation space and we normalize the highest weight state to have norm $\N w_0 \N =1$.  From~\eqref{Casimir2} we get immediately the Casimir eigenvalue:
\beq
0 = 2 \big( w_0 \,, \,  J^- J^+ w_0\big) = s(s-1) - \Omega \;\; \Rightarrow \quad
\Omega= s(s -1)\,.
\eeq
The first excited state has norm $(J^-w_0,\,J^-w_0) = (w_0,\,[J^+,J^-]w_0) = + s$
whence we conclude that unitarity requires real $s$ and 
\beq\label{s}
s>0.
\eeq

Continuing in this way, we see that the $J_3$ eigenstates $w_{-n} = (J^-)^n w_0$ in this 
representation can be labeled  by their $J_3$ eigenvalues $\big\{ |\!-s-n\rangle \big\}$ 
with $0\leq n\in\ints$. Since $J_3$ is self-adjoint, eigenvectors with distinct eigenvalues are
orthogonal. Using $(s+n)(s+n-1) - s(s-1) = n(2s-1+n)$ we find that 
\beq
\Big( (J^-)^n w_0\,,\, (J^-)^n w_0\Big) = \frac12 [n(2s-1+n)] \Big( (J^-)^{n-1} w_0\,,\, (J^-)^{n-1} w_0\Big)
\eeq
(for $n\geq 1$), so all eigenstates have positive norm squared
\beq\label{DiscNorm}
\N w_{-n} \N^2 = \prod_{k=1}^n \frac12 k(2s-1+k) = \frac{n!}{2^n} \prod_{k=1}^n (2s-1+k)\,.
\eeq
Hence the representation is unitary.
Similarly, for lowest weight representations we have $J^- w_0 = 0$ and 
$J_3 \, w_0 = +s w_0$ and the $J_3$ eigenstates $w_{n} = (J^+)^n w_0$ 
are $\big\{ |s+n\rangle\big\}$ with $0\leq n\in\ints$.  We thus see that for 
both highest and lowest weight representations unitarity is implied
by (\ref{s}). 

There are at this point no further restrictions besides $s>0$; in particular, there is 
no {\em a priori} reason to exclude {\em non-integer} values of $s>0$. It is only 
for {\em single-valued} representations of $SO(2,1)$ that one must have $\exp(2\pi i J_3)=1$, in which 
case $s$ must be a positive integer. Other \textit{rational} values can occur for covers of $SO(2,1)$ and, we will see that, in fact, all  values $s \in \rats_+$ can appear in the KM algebra.

Note that the vectors $w_{\pm n}$ are not normalized, so we adopt the notation
\beq\label{ONB}
v_n \,\coloneqq\, \frac{w_n}{\N w_n\N}
\eeq
for the rescaled and normalized vectors which give an orthonormal basis of the relevant representation space.
We also define the following notation for these discrete representation spaces:
\beq\label{DiscReps}
\cD_s^\pm = \bigoplus_{0\leq n\in\ints} \cx v_{\pm n}.
\eeq
Since all scalars are real in the formulas for the actions of $J_3$ and $J^\pm$ on these representations, we could have used $\reals v_{\pm n}$ in
the above summation and gotten a real Hilbert space of $\mf{so}(2,1)$ in this basis. Changing to the standard basis using the formulas of appendix~\ref{sec:sl2} shows, however, that we should really consider the representation space as a complex Hilbert space. The Casimir 
$\Omega = s(s-1)$ can be negative (but $\geq -1/4$) for general discrete series representations of the cover.

\subsection{Principal series representations}

Next we consider a {\it continuous} irreducible representation, $\cP_s$, labeled by a complex parameter $s$ such that the Casimir $\Omega$ acts as the scalar $s(s-1)$, 
and for which the spectrum of $J_3$ is unbounded from both above and below. 
Hermiticity of $J_3$ implies that its spectrum is real and irreducibility that all weight multiplicities are equal to one. 
Assume there exists a `minimal' eigenvector $w_p$ in $\cP_s$ satisfying $J_3\, w_p = p\, w_p$ with $0\leq p < 1$. Since we 
do not require the representations to be single-valued there does not need to be a spherical vector\footnote{For $p=0$ the vector $w_p$ is called 
a `spherical vector'.} and normalized to unity, {\em viz.}
\begin{align}
\N w_p \N^2 = 1\,.
\end{align}
Since $J_3$ is \Herm  (\ref{Herm}), $p$ is real, and so will be all other eigenvalues. 
It follows from the bracket relations (\ref{SO2,1}) that $w_{p\pm 1} = J^\pm w_p$ are 
also eigenvectors for $J_3$, but they are no longer normalized to unity.
Applying the raising and lowering operators repeatedly, we get a $J_3$ eigenbasis starting from $w_p$
\begin{align}
w_{p+n} \coloneqq \lp J^+ \rp^n w_p \;\;, \quad w_{p-n} \coloneqq \lp J^- \rp^n w_p \quad\hbox{for}\quad  n\in\nats
\end{align}
with spectrum of eigenvalues $\{p+n\mid n\in\ints\}$. 
The norms of these eigenvectors can be calculated inductively 
using $J_3\, w_p = p\, w_p$. Using (\ref{Casimir2}) 
one easily proves the following recursion relation
\begin{align}
\N w_{p+n} \N^2 &= \left( \lp J^+ \rp^n w_p ,\, \lp J^+ \rp^n w_p \right)\nn\\
&= \left( \lp J^+ \rp^{n-1} w_p ,\, (J^- J^+) \lp J^+ \rp^{n-1} w_p \right)\nn\\
&= \frac12 \big[ (p + n)(p+n-1) - s(s-1) \big] \N w_{p+n-1}\N^2
\end{align}
Similarly, we find
\beq
\N w_{p-n} \N^2 = \frac12 \big[ (p - n)(p-n-1) - s(s-1) \big]\, \N w_{p-n+1} \N^2.
\eeq
In fact, it is sufficient to just prove that for any $J_3$ eigenvector, $w_x$ with eigenvalue, $x$, 
\begin{equation}
\N w_{x\pm 1} \N^2 = \frac12 ( x(x\pm 1) - s(s-1) )\, \N w_{x} \N^2
\end{equation}
(this relation follows directly from (\ref{Casimir2})).
Therefore, for unitarity we must have $(p + n)(p + n - 1) > s(s-1)$ for all $n\in\ints$, 
so if $x = p+n$, we may choose $n$ such that
$1 > x \geq 0$ so $0 > x-1 \geq -1$ so $0 \geq x(x-1) > s(s-1)$. 
The minimum value of the parabola $y = x(x-1)$ is $-1/4$, so we will be sure of this condition for all $n$ and $p$ 
when $\Omega = s(s-1) < -1/4$. Then all norms are positive and we obtain a unitary representation. Below we write $\Omega$ for its value on $\cP_s$. 

Among the continuous representations one further distinguishes 
between {\bf principal series} for which $\Omega\leq -\frac14$, and 
{\bf complementary series} for which $-\frac14 < \Omega < 0$. The distinction arises since 
$-1/4$ is the minimum of the parabola $s(s-1)$ for real $s$; the unitary principal series requires complex $s$. 
The relation above gives the following closed formula for the norms of the eigenvectors: 
\begin{align}\label{Norms}
\N w_{p\pm n}\N^2 = \frac{1}{2^{n}} \prod_{k=0}^{n-1} \big[(p\pm k)(p\pm (k+1) - \Omega\big] 
\end{align}
Writing $\Omega=s(s-1)$ also for the principal series representations we have $s=\frac12(1+ iq)$ with $q\in\reals$, thus $\Omega=-\frac14(1+q^2) = -s{\bar s}$. 
The parameter $q$ is determined up to sign and its value depends on the example. 

Defining the orthonormal basis as in (\ref{ONB}) we denote denote the corresponding
representation space by 
\beq\label{PrincReps}
\cP_s = \bigoplus_{n\in\ints} \cx v_{p+n}\,.
\eeq

\subsection{Finite-dimensional representations}

For completeness we also mention finite-dimensional irreducible modules of $\mf{sl}(2,\mathbb{R})$, see for instance~\cite{Humphreys}. These will be denoted by $V(m)$ for the $(m+1)$-dimensional module and the only unitary case is the trivial representation $V(0)$. The non-unitary ones play a role in the decomposition with respect to real roots $\alpha$, discussed for example in section~\ref{sec:Fibreal}.

The module $V(m)$ can be characterised by having a highest weight vector $v_m$ with eigenvalue $m/2$ under $J_3$ and satisfying $J^+v_m=0$. The Casimir eigenvalue on such an irreducible representation is $\Omega= \tfrac14 m(m+2)$ in the normalization~\eqref{Casimir1}. In particular, for the adjoint representation with $m=2$, the Casimir is $\Omega=2$. The Casimir spectrum on the non-unitary representations $V(m)$ overlaps with that of the unitary representations.

\section{Hilbert space realizations}
\label{app:FR}

For convenience we here summarize the known Hilbert (function) space realizations
of the unitary representations of SL(2,$\reals$) and of its covers.

\subsection{\texorpdfstring{Functional representations of $SL(2,\mathbb{R})$}{Functional representations of SL(2,R)}}

For all the discrete series representations of the group SL(2,$\reals$) on Hilbert spaces of complex-valued square integrable functions, $G(z)$, 
the left action is defined by
\beq\label{Mobiusf}
G(z) \;\; \rightarrow \;\; (S\cdot G)(z) \coloneqq G\left(\frac{az+b}{cz+d}\right) (cz+d)^{-2s} \;\;\;,\quad\hbox{for}\quad
S^{-1} =\begin{bmatrix}
a & b\\
c & d
\end{bmatrix} \in {\rm SL}(2,\reals)
\eeq
with \textit{integer} $s$.
For the continuous series representations the factor $(cz+d)^{-2s}$ is replaced by $|cz+d|^{-2s}$. 
With $S_t$ taken from each of the one-parameter subgroups, $\{\exp(te)\mid t\in\reals\}$, $\{\exp(tf)\mid t\in\reals\}$ and $\{\exp(th)\mid t\in\reals\}$, 
for fixed $z$, the linear term in the Taylor expansion of $L(t) = (S_t\cdot G)(z)$ obtained from $\lim_{t\to 0}L'(t)$, gives the differential operators 
representing the $\mf{sl}(2,\reals)$ basis vectors (\ref{efhbasis}) to be
\beq\label{efhops}
E = -\frac{d}{dz}, \qquad F = z^2 \frac{d}{dz} + 2sz, \qquad H = -2z\frac{d}{dz} - 2s
\eeq
where we use capital letters to distinguish these operators from the abstract Lie algebra elements they represent. 
The basis vectors\footnote{This change of basis is similar to the change to the so-called `compact basis' that appears for example in~\cite{Howe:1992}. However, there are a few sign differences and a rescaling by a factor $\sqrt{2}$ to obtain~\eqref{SO2,1}.}
\beq\label{conjbasis}
j_3 = \frac{\bi}{2}(e - f), \qquad j^+ = \frac{1}{2\sqrt{2}} (-\bi(e + f) + h), \qquad j^- = \frac{1}{2\sqrt{2}} (-\bi(e + f) - h)
\eeq
in $\mf{sl}(2,\cx)$ satisfy the (\ref{SO2,1}) Lie brackets $[j_3, j^\pm] = \pm j^\pm$ and $[j^+, j^-] = -j_3$ but also satisfy $X^\dagger\eta + \eta X = 0$ for 
$\eta = \left[\begin{smallmatrix} 0&1\\1&0\end{smallmatrix}\right]$, which means this real Lie algebra is $\mf{su}(1,1)$. The differential operators corresponding to these complex
linear combinations of the operators in (\ref{efhops}) give us an explicit  
realization of the algebra (\ref{SO2,1}) by the differential operators
\begin{align}
\label{eq:dopnew}
J_3 &\,=\, 
    -\frac{i}{2}(1+z^2) \frac{d}{dz} -isz\,,\nn\\
J^+ &\,=\,
\frac{1}{2\sqrt{2}} \lp -i (z-i)^2\frac{d}{dz} -2s(1+iz)\rp \,,\nn\\
J^- &\,=\,
\frac{1}{2\sqrt{2}} \lp  -i (z+i)^2\frac{d}{dz} +2s(1-iz)\rp .
\end{align}
One still has to specify the Hilbert space on which these operators act, and the scalar product
$(\cdot ,\cdot)$ that defines the norm. The variable $z$ is real for the continuous representations
and complex for the discrete representations (in which case we write $z = x +iy$).

For all realizations we insist on the hermiticity properties  
\begin{align}\label{Herm1}
(J_3)^\dagger \,=\, J_3 \,,\quad
(J^\pm)^\dagger \,=\, J^\mp
\end{align}
where hermiticity is defined with respect to the given scalar product. From these expressions it is 
straighforward to compute the Casimir
\begin{align}
\Omega &\,=\,  -J^+J^- -J^-J^+ + J_3 J_3 = s(s-1)
\end{align}

In terms of the operators (\ref{SO2,1}) this action corresponds to the exponential
\beq
\exp \big(iw J^+ + i\bar{w} J^- + ir J_3\big)
\eeq
with appropriate parameters $w = u + iv\in\cx$ and $r \in\reals$
determined from the matrix $S$. 
Note that if we use the formulas (\ref{conjbasis}), the expression
\beq\label{realexpression}
iw j^+ + i\bar{w} j^- + ir j_3 = \frac{1}{\sqrt{2}} (u(e + f) - vh) +  \frac{1}{2} r (e - f)
\eeq
is a real linear combination of $\{e,f,h\}$, so it is in $\mf{sl}(2,\reals)$ and its 
exponential is in $SL(2,\reals)$. 
We note that with respect to the Kac bilinear form we have that $H^\dagger= - H$, $E^\dagger=-E$ and $F^\dagger=-F$.

The map $S\in SL(2,\reals)$ is represented with respect to an orthonormal basis of functions, $v_n$, 
by a unitary matrix
\beq\label{S}
S \cdot v_n = \sum_{m\in\ints}  {\rm U}_{mn}(S) v_m\,.
\eeq
Because the Hilbert space consists of complex valued functions the infinite
matrix U$_{mn}$ is necessarily complex. Unitarity means that
\begin{align}
\sum_{k\in \ints} {\rm U}^*_{km} {\rm U}_{kn} = \delta_{mn}
\end{align}
so in particular all sums that arise are manifestly convergent.

\subsubsection{Discrete series}

As explained above, the function space realizations of discrete series representations correspond to lowest or highest weight representations,
depending on whether there is a ground state killed by $J^-$ or by $J^+$, respectively. In the first case, $z$ is a {\em complex} variable 
$z = x+iy$ and the functions $f(z), g(z)$ are holomorphic in the  upper  half-plane $\UHP$.
The scalar product is
\beq\label{dz}
(f,g) \coloneqq \int_{\UHP} dx dy \, y^{2s-2} \overline{f(z)} g(z)
\eeq
Importantly, the operators (\ref{eq:dopnew}) are \Herm  with respect to (\ref{dz}) only for $s>1$ because
only then the boundary term arising from  integration by parts vanishes.

The normalized ground state $\vp_0$ of a lowest weight representation satisfying  $J^- \vp_0= 0$ and $J_3\vp_0=s\vp_0$ is given by
\beq
\vp_0(z) = \sqrt{\frac{2s-1}{\pi}} \ 2^{2s-1}\ \frac1{(z+i)^{2s}} 
\eeq
because 
\beq
 \int_{x\in\reals} \int_{y>0} dxdy\ y^{2s-2}\ (x^2 + (y+1)^2)^{-2s} = \frac{4^{1-2s}\pi}{2s-1}
 \eeq
for $s > \frac12$, and the excited state $\vp_n$ with $J_3$-eigenvalue $n+s$ is
\beq\label{phin}
\vp_n(z) = (J^+)^n \vp_0(z) \,=\, A_n \frac{(z-i)^n}{(z+i)^{2s + n}} 
\eeq
with 
\beq
\label{eq:anconst}
A_0 =  \sqrt{\frac{2s-1}{\pi}} \ 2^{2s-1}\quad\hbox{and}\quad A_n = \frac{A_0}{\sqrt{2^n}} \prod_{k=0}^{n-1} (2s+k) = \frac{A_0}{\sqrt{2^n}} \frac{\Gamma(2s+n)}{\Gamma(2s)}
\eeq
because we have the recursion 
\beq
\label{eq:normd}
\N \vp_n \N^2 = \frac12 n(2s-1+n) \N \vp_{n-1} \N^2 
= 2^{-n} \frac{\Gamma(n+1)\Gamma(2s+n)}{\Gamma(2s)} \N \vp_0 \N^2\,.
\eeq
The functions $\vp_n$ are not normalized. 

Similar definitions apply for highest weight representations with a ground state such that $J^+ \vp_0= 0$ and $J_3\vp_0 = -s\vp_0$, where the integral 
(\ref{dz}) is now to be performed over the {\em lower} half-plane, and we write $\vp_{-n} = (J^-)^n \vp_0(z)$ for $n\in\nats$. 

We use the following notation for these two kinds of function spaces which are each irreducible $\so$ representation spaces
\beq
D_s^\pm = \bigoplus_{n\geq 0} \cx \vp_{\pm n}(z)
\eeq
but $\cx$ could be replaced by $\reals$ since all the coefficients $A_n$ are real. We wish to define a Lie algebra module 
isomorphism which is also a (real or complex) Hilbert space isometry
\beq
\Phi_s^\pm : \cD_s^\pm \to D_s^\pm.
\eeq
Using the orthonormal basis of $\cD_s^\pm$ in \eqref{DiscReps}, 
define $\Phi_s^\pm(v_0) = \vp_0(z)$ and require that 
\beq
\Phi_s^\pm((J^\pm )^n v_0) = (J^\pm )^n \Phi_s^\pm( v_0) = (J^\pm )^n \vp_0(z)\,.
\eeq

For completeness let us mention that the functional realization of the discrete series
representation can be equivalently done on the Poincar\'e disk $|w| <1$ by means of the 
standard Cayley transformation
\begin{align}
w = \frac{z-i}{z+i} 
\quad \Longleftrightarrow \quad
z = i \frac{1+w}{1-w}
\end{align}
mapping the upper half plane (with coordinate $z$) to the interior of the unit disk (with coordinate $w$).
The differential operators~\eqref{eq:dopnew} can be easily converted using the Cayley transformation, leading for example to
\begin{align}
J_3 = w \frac{d}{dw} + s \frac{1+w}{1-w}\,.
\end{align} 
The integral~\eqref{dz} becomes for $w=u+i v$ and $s$ a positive integer
\begin{align}
\label{eq:ipPD}
(f,g) =  \int_{|w|<1} du dv \, 
\frac{(1-|w|)^{2s-2}}{(1-w)^{2s}(1-\bar{w})^{2s}} \overline{f(w)} g(w)
 \,.
\end{align}
The $J_3$-eigenfunctions $\vp_n$ in (\ref{phin})
are then replaced by functions proportional to $(1-w)^{2s} w^n$. 
In the inner product~\eqref{eq:ipPD} such functions have the norm
\begin{align}
 \int_{|w|< 1} du dv \, (1 - |w|^2)^{2s-2} |w|^{2n} = \frac{\pi}{2s-1}  \frac{\Gamma(2s)\Gamma(n+1)}{\Gamma(2s+n)}
\end{align}
from a standard representation of the Euler beta function.
When discussing covers of $SO(2,1)$ this formula must be analytically continued to non-integer values of $s>0$
and corresponding Hilbert spaces can be defined \cite{Sally:1965},
see also section~\ref{app:covers}. We also note that by omitting the factor
$(2s-1)^{-1}$ one can redefine the norm to be
finite and positive for all $s>0$ \cite{Sally:1965}.

\subsubsection{Principal series representations}

A concrete realization of the principal series representations of $SL(2,\mathbb{R})$  is provided by the
space L$^2(\reals)$ of complex valued square integrable functions over the real line 
$\reals$ (so $x$ is now real) with the scalar product
\begin{align}
\label{eq:ipPS}
(f,g) \,\coloneqq \, \int_{-\infty}^\infty \overline{f(x)} g(x) dx
\end{align}
with the above operator realizations. 
Keeping in mind that $( (x+a)^2 \frac{d}{dx})^\dagger = -(x+\bar{a})^2 \frac{d}{dx} -2(x+\bar{a})$, 
with an extra linear term, one sees that the hermiticity properties (\ref{Herm1}) are only satisfied  if
\beq
s=\frac12(1+iq) \;\; \Rightarrow \quad \Omega = s(s-1) = -\frac14(1 + q^2) \leq - \frac14
\eeq
with $q\in \reals$.
The $J_3$-eigenfunctions are given by
\beq\label{J3efunctions}
J_3 \vp_m(x) = m \vp_m(x) \;\; \Rightarrow \quad
\vp_m(x) = \frac1{\sqrt{\pi}} \frac1{(1 + x^2)^s} \ \left[\frac{i-x}{i+x}\right]^m\,.
\eeq
Note that we are here considering representations of $SL(2,\mathbb{R})$ so that $m\in\mathbb{Z}$.
Furthermore, $\N \vp_m(x) \N = 1$ because
\begin{equation}
\int_{-\infty}^\infty \overline{\vp_m(x)} \vp_m(x) dx 
= \frac{1}{\pi}\ \int_{-\infty}^\infty \frac{dx}{1 + x^2} = 1.
\end{equation}
To solve the differential equation $J_3 \vp_m(x) = m \vp_m(x)$ we have also used the identity
\begin{equation}
\exp\big(2i\arctan(x)\big) = \frac{i-x}{i+x}.
\end{equation}
When we use the differential operators representing $J^\pm$ from (\ref{eq:dopnew}), we get 
\begin{equation}
J^\pm \vp_m(x) = -\frac{(m \pm s)}{\sqrt{2}} \vp_{m\pm 1}(x)
\end{equation}
It follows that for $n\geq 1$ we have
\beq
(J^+)^n \vp_m(x) = A_{n,s} \vp_{m+n}(x)  \quad\hbox{and}\quad  
(J^-)^n \vp_m(x) = A_{-n,s}  \vp_{m-n}(x) 
\eeq
where
\beq
A_{\pm n,s} = (-1)^n \frac{(m\pm s)(m\pm 1\pm s)\cdots(m\pm (n-1)\pm s)}{(\sqrt{2})^n}
\eeq
For $0\leq p <1$ denote by $\Pi_s^{(p)} = \bigoplus_{n\in\ints} \cx \vp_{p+n}$ 
the irreducible representation space of functions (with $m$ replaced by $p+n$ in \eqref{J3efunctions}) we have just found for the Lie algebra 
$\so$ with operator basis given in (\ref{eq:dopnew}). Using the basis $\{w_{p+n}\mid n\in\ints\}$ of $\cP_s$ in \eqref{Norms}, we can define a module isomorphism 
$\Phi_s: \cP_s \to \Pi_s^{(p)}$ which is also a Hilbert space isometry, as follows. Begin by setting 
$\Phi_s(w_p) = \vp_p$. Then let
\begin{align}
\Phi_s(w_{p\pm 1}) &= \Phi_s(J^\pm w_{p}) = J^\pm  \Phi_s(w_{p}) = 
J^\pm \vp_p = - \frac{(p \pm s)}{\sqrt{2}} \vp_{p\pm 1},\\
\Phi_s(w_{p\pm 2}) &= \Phi_s(J^\pm w_{p\pm 1}) = J^\pm  \Phi_s(w_{p\pm 1}) = \frac{(p \pm s)}{\sqrt{2}} J^\pm \vp_{p\pm 1} =  
\frac{(p \pm s)}{\sqrt{2}} \frac{(p\pm 1 \pm s)}{\sqrt{2}}\vp_{p\pm 2}\nn
\end{align}
so in general, for any $n\in\nats$, let
\begin{equation}
\Phi_s(w_{p\pm n}) = \Phi_s((J^\pm)^n w_{p}) = (J^\pm)^n  \Phi_s(w_{p}) = (J^\pm)^n \vp_p = A_{\pm n,s} \vp_{p\pm n}.
\end{equation}
To check that this is a Hilbert space isometry we need to see that 
\begin{equation}
\N w_{p\pm n} \N^2 = \N A_{\pm n,s} \vp_{p\pm n} \N^2 = A_{\pm n,s}\ {\overline{A_{\pm n,s}}} = 
\frac{1}{2^n}\ \prod_{k=0}^{n-1} (p\pm k\pm s) (p\pm k\pm {\bar s})
\end{equation}
matches the answer in (\ref{Norms}). For each $k$ we have 
\begin{equation}
(p\pm k\pm s) (p\pm k\pm {\bar s}) = (p\pm k)^2 \pm (p\pm k)(s + {\bar s}) + s {\bar s} = (p\pm k)^2 \pm (p\pm k) + s {\bar s} 
\end{equation}
does match
\begin{equation}
((p\pm k)(p\pm (k+1) - \Omega) = ((p\pm k)(p\pm k \pm 1) + s {\bar s}) .
\end{equation}

\subsubsection{Complementary series}

Finally we turn to the complementary series with real $0<s<1$.  In this case we have Casimir eigenvalues $-\tfrac14<\Omega<0$.
 These are
discussed (in the compact $SU(1,1)$ picture) in~\cite{Sugiura:1990}, 
see also~\cite{Knapp:1986,Lang:1985,Howe:1992} for further discussions.
Even though we have not found any Kac--Moody algebra where they arise, we give some details here for completeness.

A model is provided by complex-valued functions on the real line.
The inner product is given by
\begin{align}
(f,g) = \int_{\reals^2} \overline{f(x)} g(y) |x-y|^{2s-2} dx dy
\end{align}
and the Hilbert space is the completion with respect to this norm. We will show this
in the appendix assuming $f,g \in$ L$^1(\reals) \cap$ L$^2(\reals)$. We will now show
that we have
\begin{itemize}
\item convergence
\item positive definite
\item operators have right hermiticity properties
\end{itemize}
provided
\beq
0<s <1 \;\; \Rightarrow \quad    -\frac14 < \Omega < 0
\eeq
Because $f$ and $g$ are in L$^1(\reals) \, \cap$ L$^2(\reals)$ we can write them as Fourier integrals
\begin{align}
f(x) = \int_\reals a(k) e^{2\pi i k x} dk \,,\quad
g(x) =\int_\reals b(k) e^{2\pi i k x} dk \,.
\end{align}
Then, after changing variables
\begin{align}
v = \frac12(x-y)\,,\quad u = \frac12(x+y)\,,
\end{align}
one can rewrite the integral
\begin{align}\label{fg}
(f,g)_\alpha &= \int_{\reals^2} \overline{f(x)} g(y) |x-y|^\alpha dx dy\nn\\
&=  \int_{\reals^2} \int_{\reals^2} \overline{a(k)} b(\ell) e^{2\pi i( kx -\ell y)} |x-y|^\alpha dx dy dk d\ell\nn\\
&=  2^{\alpha+1} \int_{\reals^2} \int_{\reals^2} \overline{a(k)} b(\ell) e^{2\pi i( k (u+v) -\ell(u-v))} |v|^\alpha 
du dv dk d\ell\nn\\
&=   \int_\reals \int_{\reals} \overline{a(k)} b(k) e^{2\pi i kv} |v|^\alpha dv dk\nn\\
&=2^{-\alpha} \pi^{-1-\alpha}\Gamma(1+\alpha)\sin\lp \frac{\pi|\alpha|}{2}\rp \int |k|^{-1-\alpha} 
\overline{a(k)} b(k) dk\,,
\end{align}
where we used  
\begin{align}
\int_{\reals}  e^{2\pi i nv} |v|^\alpha dv = -2^{-\alpha} |n|^{-1-\alpha} \pi^{-1-\alpha} \Gamma(1+\alpha) \sin\lp \frac{\pi\alpha}{2}\rp\,.
\end{align}
For convergence at small $k$ and large $v$ we must have
\beq 
-1<\Re(\alpha)<0
\eeq
where the lower bound comes from the pole of the $\Gamma$-function at 0, while the 
upper limit comes from requiring convergence of the integral over Fourier coefficients
in (\ref{fg}). In particular, it follows from these formulas that
$(f,f)_\alpha$ converges for this choice and is positive definite. 

The admissible values for $\alpha$ in terms of $s$ can be fixed by demanding that the 
operators be Hermitian. To this aim consider the adjunction of the operators~\eqref{eq:dopnew} 
with respect to this inner product, assuming that $s\in\reals$. For $J_3$ one gets:
\begin{align}
(f, J_3 g)_\alpha &=  \int_{\reals^2} \overline{f(x)} \lp -\frac{i}2 (1+y^2) \frac{d}{dy} g(y) - i s y g(y)\rp |x-y|^\alpha dx dy\nn\\
& =  \int_{\reals^2} \overline{f(x)} g(y) \lp \frac{i}2 (1+y^2) \frac{d}{dy} |x-y|^\alpha + i y (1-s) |x-y|^\alpha \rp dx dy\,.
\end{align}
Compute similarly, using $\frac{d}{dx} |x-y|^\alpha = -\frac{d}{dy}|x-y|^\alpha$, that
\begin{align}
(J_3 f, g)_\alpha &=  \int_{\reals^2} \lp \frac{i}2 (1+x^2) \frac{d}{dx} \overline{f(x)}  + i s x \overline{f(x)}\rp  g(y) |x-y|^\alpha dx dy\nn\\
& =  \int_{\reals^2} \overline{f(x)} g(y) \lp \frac{i}2 (1+x^2) \frac{d}{dy} |x-y|^\alpha - i x (1-s) |x-y|^\alpha \rp dx dy\,.
\end{align}
Therefore the difference is
\begin{align}
(f, J_3 g)_\alpha - (J_3f,  g)_\alpha &= \int_{\reals^2} \overline{f(x)} g(y) \lp \frac{i}{2} (y^2-x^2) \frac{d}{dy} |x-y|^\alpha +i(1-s)(x+y) |x-y|^\alpha \rp dx dy\,.
\end{align}
Let us focus on the term in parentheses:
\begin{align}
&\quad\quad \frac{i}{2} (y^2-x^2) \frac{d}{dy} |x-y|^\alpha +i(1-s)(x+y) |x-y|^\alpha\nn\\
&=\frac{i}{2} \alpha (x-y) (x+y) |x-y|^{\alpha-1} \sgn(x-y) +i (1-s) (x+y) |x-y|^\alpha\nn\\
&=\lp \frac{i}{2} \alpha +i (1-s) \rp (x+y) |x-y|^{\alpha} \,,
\end{align}
such that this vanishes when 
\begin{align}
\label{eq:alphas}
\alpha = 2s-2\,,
\end{align}
making $J_3$ \Herm  with respect to this inner product.

The condition $-1<\Re(\alpha)<0$ translates for real $s$ into
\begin{align}
\frac12<s<1\,,
\end{align}
which is half of the interval $0<s<1$ and the other half is covered by analytic continuation using the functional relation between $s$ and $1-s$.

We also demonstrate the other Hermiticity relation. For $J^+$ one finds:
\begin{align}
(f, J^+ g)_\alpha &=  \frac{1}{2\sqrt{2}} \int_{\reals^2} \overline{f(x)} \lp -i(y-i)^2 \frac{d}{dy} g(y) - 2 s (1+i y) g(y)\rp |x-y|^\alpha dx dy\nn\\
& =  \frac{1}{2\sqrt{2}} \int_{\reals^2} \overline{f(x)} g(y) \lp i (y-i)^2 \frac{d}{dy} |x-y|^\alpha -2 i (s-1) (y-i)|x-y|^\alpha \rp dx dy\,.
\end{align}
and
\begin{align}
(J^- f, g)_\alpha &=  \frac{1}{2\sqrt{2}} \int_{\reals^2} \lp i (x-i)^2 \frac{d}{dx} \overline{f(x)}  + 2 s (1+ix) \overline{f(x)}\rp  g(y) |x-y|^\alpha dx dy\nn\\
& = \frac{1}{2\sqrt{2}} \int_{\reals^2} \overline{f(x)} g(y) \lp i (x-i)^2 \frac{d}{dy} |x-y|^\alpha + 2i (s-1)(x-i) |x-y|^\alpha \rp dx dy\,.
\end{align}
Therefore the difference is
\begin{align}
(f, J^+ g)_\alpha - (J^-f,  g)_\alpha &= \frac{1}{2\sqrt{2}} \int_{\reals^2} \overline{f(x)} g(y) \bigg( -i(x-y)(x+y-2i)\frac{d}{dy} |x-y|^\alpha\\
&\hspace{40mm} -2i(s-1)(x+y-2i) |x-y|^\alpha \bigg) dx dy\,.\nn
\end{align}
The term in parentheses becomes
\begin{align}
&\quad -i(x-y)(x+y-2i)\frac{d}{dy} |x-y|^\alpha-2i(s-1)(x+y-2i) |x-y|^\alpha \nn\\
&= i\alpha (x-y)(x+y-2i) |x-y|^{\alpha-1} \sgn(x-y)-2i(s-1)(x+y-2i) |x-y|^\alpha \nn\\
&= i(x+y-2i) |x-y|^\alpha ( \alpha -2(s-1)) \,,
\end{align}
so that for $\alpha=2s-2$ one has $(J^+)^\dagger = J^-$ as required and consistent with~\eqref{eq:alphas}.

\subsection{Representations of covers}
\label{app:covers}

The functional realization on covers proceeds by a very similar method. We explain this in the case of the discrete series on an $N$-fold cover $\tilde{G}$ of $SL(2,\mathbb{R})$, meaning a $2N$-fold cover of $SO(2,1)$. Having an $N$-fold cover, means that the group $\tilde{G}$ formally consists of pairs $(S,\zeta)$ with $S\in SL(2,\mathbb{R})$ and $\zeta$ is an element of the cyclic group of order $N$ that can be represented by $N$th roots of unity. The product on such pairs involves a cocycle that defines the extension. We will not require the precise form of this cocycle (see e.g.~\cite{Hauser:2023}) since its definition will be implicit in our construction.

Let 
\begin{align}
S^{-1} = \begin{bmatrix}a & b\\c &d\end{bmatrix} \in SL(2,\mathbb{R})
\end{align}
and choose an $N$th root $\mu$ of the linear function $cz+d$ where $z$ is in the upper half plane. We demand that $\mu$ is holomorphic on the upper half plane, and by construction $\mu(z)^N = cz+d$. Such a function is given up the choice of a root of unity and there are $N$ such roots. Therefore the pairs $(S,\mu)$ are what is needed to describe the $N$-fold cover of $SL(2,\mathbb{R})$. 

The action of a such a pair $(S,\mu)$ on holomorphic functions $G(z)$ on the upper half plane is
\begin{align}
G(z) \mapsto ( (S,\mu)\cdot G)(z) = \frac{1}{\mu(z)^{2s}} G\left(\frac{az+b}{cz+d}\right)
\end{align}
and generalizes~\eqref{Mobiusf} to the case when $s\in \frac{1}{2N} \mathbb{N}$. 
The infinitesimal action of the Lie algebra $\mf{sl}(2,\mathbb{R})$ in terms of differential operators is unchanged from this definition, in agreement with the fact that the differential operators~\eqref{efhops} satisfy the $\mf{sl}(2,\mathbb{R})$ Lie algebra for any complex $s$.

From this action one deduces the product on the covering group
\begin{align}
(S_1,\mu_1) \cdot (S_2, \mu_2) = (S_1S_2, (\mu_1\circ S_2) \mu_2)\,,
\end{align}
where the second entry denotes the holomorphic function
\begin{align}
\left((\mu_1\circ S_2) \mu_2 \right)(z) = \mu_1\left( \frac{a_2z+b_2}{c_2z+d_2}\right) \mu_2(z)
\end{align}
that intertwines the product of the roots with the action of $SL(2,\mathbb{R})$ on the upper half plane. 

The Hilbert space of the discrete series of $\tilde{G}$ consists of all holomorphic functions with finite norm with respect to the (analytically continued) norm~\eqref{dz}. For fractional $0<s<1$ this represents a more stringent requirement than holomorphicity on the upper half but has been discussed in detail in the literature, see for example~\cite{Sally:1965}. In particular, the hermiticity properties~\eqref{Herm} and unitarity of the representation can be maintained.

For the principal series the holomorphic functions on the upper half plane are taken to the boundary real line. 

\newpage

\section{Figures}

\newpage

\newpage

\newpage

\newpage

\begin{thebibliography}{56}
\bibitem{Tsu} H.~Tsurusaki,
``$\mf{sl}_2$ triples whose nilpositive elements are in a space which is spanned 
by the real root vectors in rank-2 symmetric hyperbolic Kac-Moody Lie algebras,''
\eprintRT{2107.05234}.
\end{thebibliography}
\end{document}